\documentclass{siamart1116}
\usepackage{threeparttable}
\usepackage{lipsum}
\usepackage{amsfonts}
\usepackage{graphicx}
\usepackage{subfigure}
\usepackage{diagbox}
\usepackage{epstopdf}
\usepackage{algorithmic}
\usepackage{amsfonts,amssymb,mathrsfs,amsmath,color,fancyhdr}
\usepackage{multirow}
\usepackage{psfrag}
\usepackage{graphics}
\usepackage{float}
\usepackage{multirow}
\usepackage{enumerate}
\usepackage{caption}
\usepackage{booktabs}
\usepackage{array}
\usepackage{subeqnarray}
\usepackage{cases}
\ifpdf
  \DeclareGraphicsExtensions{.eps,.pdf,.png,.jpg}
\else
  \DeclareGraphicsExtensions{.eps}
\fi

\numberwithin{equation}{section}

\newcommand{\E}{\mathbb{E}}
\newcommand{\msE}{\mathsf{E}}

\newcommand{\R}{\mathbb{R}}

\newcommand{\F}{\mathcal{F}}
\newcommand{\mtF}{\mathbb{F}}

\newtheorem{sch}{ {\bf Scheme} }[section]

\newtheorem{lem}{Lemma}[section]

\newtheorem{thm}{Theorem}[section]
\newtheorem{pro}{Proposition}[section]

\newtheorem{rem}{Remark}[section]
\newtheorem{coro}{Corollary}[section]

\newtheorem{ex}{Example}[section]

\headers{Numerical methods for  MSDEJs}{Yabing Sun and Weidong Zhao}

\title{Numerical methods for  mean-field
stochastic differential equations with jumps
\thanks{This work is partially supported by the science challenge Project (No.
TZ2018001), the NSF of China (under Grant Nos. 11571351, 11571206, 11831010, 11871068)
and the China Postdoctoral Science Foundation (No. 2019TQ0073).}}
\author{Yabing Sun\footnotemark[1]
\and Weidong Zhao\footnotemark[2]
}

\begin{document}
\maketitle

\renewcommand{\thefootnote}{\fnsymbol{footnote}}
\footnotetext[2]{Corresponding author}
\footnotetext[1]{College of Science, National University of Defense Technology,
Changsha, Hunan 410073, China. Email: sunybly@163.com}
\footnotetext[2]{School of Mathematics \& Institute of Finance,
Shandong University, Jinan, Shandong 250100, China.
 Email: wdzhao@sdu.edu.cn}

\begin{abstract}
In this paper, we are devoted to the numerical methods
for mean-field stochastic differential equations with jumps (MSDEJs).
First by using the mean-field It\^o formula
[Sun, Yang and Zhao, Numer. Math. Theor. Meth. Appl., 10 (2017), pp.~798--828],
we  develop the It\^o formula and construct the It\^o-Taylor expansion for MSDEJs.
Then based on the It\^o-Taylor expansion,
we propose the strong order $\gamma$ and the weak order $\eta$
It\^o-Taylor schemes for MSDEJs.
The strong and weak convergence rates $\gamma$ and $\eta$
of the strong and weak It\^o-Taylor schemes are theoretically proved, respectively.
Finally some numerical tests are also presented to verify our theoretical conclusions.
\end{abstract}



\section{Introduction}

Let $(\Omega,\mathcal{F},\mathbb{F},P)$ be a complete filtered probability space
with $\mathbb{F}=\{\mathcal{F}_t\}_{0\le t\le T}$ being the filtration of
the following two mutually independent stochastic processes:
\begin{itemize}
\item the $m$-dimensional Brownian motion: $W = (W_t)_{0\leq t\leq T}$;
\item the Poisson random measure on $\mathsf{E}\times[0, T]$:
\{$\mu(A\times[0,t])$, $A\in\mathcal{E}$, $0\le t\le T$\},
where $\mathsf{E} = \R^q\backslash\{0\}$ and  $\mathcal{E}$ is its Borel field.
\end{itemize}
Suppose that $\mu$ has the intensity measure
$\nu(de,dt)=\lambda(de)dt$, where $\lambda$ is a $\sigma$-finite measure on $(\mathsf{E},\mathcal{E})$
satisfying $\int_{\mathsf{E}} (1 \land |e|^2) \lambda (de) < + \infty$.
Then we have the compensated Poisson random measure
$$\tilde \mu(de,dt)=\mu(de,dt)-\lambda(de)dt,$$
such that $\{\tilde \mu(A\times[0,t])=(\mu-\nu)(A\times[0,t])\}_{0\le t\le T}$
is a martingale for any $A\in \mathcal{E}$ with $\lambda(A)<\infty$.
Moreover, let $F$
be the distribution
of the jump size,
then it holds that
$$\lambda(de)=\dot{\lambda} F(de),$$
where $\dot{\lambda}=\lambda(\mathsf{E})<\infty$
is the intensity of the Poisson process
$N_t=\mu(\mathsf{E}\times[0,t])$, which counts the number of jumps of $\mu$
occurring in $[0,t]$.
Then, the Poisson measure $\mu$
generates a sequence of pairs $\{(\tau_i,Y_i),i=1,2,\dots,N_T\}$
with $\{\tau_i\in[0,T],i=1,2,\dots,N_T\}$ representing the jump times of
the Poisson process $N_t$ and
$\{Y_i\in\msE,i=1,2,\dots,N_T\}$ the corresponding jump sizes
satisfying $Y_i\overset{iid}\sim F$.
For more details of the Poisson random measure
or L\'evy measure, the readers are referred to \cite{CT2004,PB2010}.

We consider the following mean-field
stochastic differential equation with jumps (MSDEJs) on
$(\Omega,\mathcal{F},\mathbb{F},P)$
\begin{equation}
\begin{aligned}\label{MSDEJ}
X_t^{t_0,\xi} = \xi &+\! \int_{t_0}^t\E\big[b(s,X_s^{t_0,\xi'},x)\big]\big|_{x=X_s^{t_0,\xi}} ds
+\! \int_{t_0}^t\E\big[\sigma(s,X_s^{t_0,\xi'},x)\big]\big|_{x=X_s^{t_0,\xi}}  dW_s \\
& +\! \int_{t_0}^t\int_{\msE} \E\big[c(s,X_{s-}^{t_0,\xi'},x,e)\big]\big|_{x=X_{s-}^{t_0,\xi}}
\mu(de,ds),\;\;\;\;\;\;\;\;\;0\leq t_0\le t\le T,
\end{aligned}
\end{equation}
where $t_0$ and $T$ are, respectively, the deterministic initial and terminal time;
the initial condition $\xi$ is $\mathcal{F}_{t_0}$ measurable;
$b:[0,T]\times\R^d\times\R^d\rightarrow\R^d$,
$\sigma:[0,T]\times\R^d\times\R^d\rightarrow\R^{d\times m}$,
and $c:[0,T]\times\R^d\times\R^d\times\msE\rightarrow\R^d$
are the so called drift, diffusion and jump coefficients, respectively.
Here the superscript $^{t_0,\xi}$
indicates that the MSDEJ \eqref{MSDEJ}
starts from the time-space point $(t_0,\xi)$,
and $X_t^{t_0,\xi'}$ is the solution
of the MSDEJ \eqref{MSDEJ} with $\xi=\xi'$.
In general, $\xi$ and $\xi'$ are different.

Mean-field stochastic differential equations (MSDEs),
also called McKean-Vlasov SDEs, was first studied by Kac \cite{Kac56,Kac58} in the 1950s.
Since then, MSDEs have been encountered and intensively investigated in many areas
such as kinetic gas theory \cite{bt,m,TV}, quantum mechanics \cite{mav}, quantum chemistry \cite{sss},
McKean-Vlasov type partial differential equations (PDEs) \cite{blp,blp2017,k,M1966},
mean-field games \cite{cf,CD2019,gll,ll}
and mean-field backward stochastic differential equations (MBSDEs) \cite{bdl,blp,BY2015,SZ2019,SZZ2018}.
%
In the last decade, MSDEJs have also received much attention
because of its wild applications in
the research on nonlocal PDEs \cite{HL2016}, MBSDEs with jumps \cite{Li2017,LM2016},
economics and finance \cite{HTB2015},
and mean-field control and mean-field games with jumps \cite{HAA2015,zjf,zjf1}.
Therefore, it is important and necessary to  study the numerical solutions
of MSDEJs.

%
Compared with the
well developed theory of numerical methods for  stochastic differential equations with jumps (SDEJs)
(see \cite{HK2005,LL2010,PB2010} and references therein),
little attention has been paid to the numerical methods for MSDEJs.
In this work, we aim to propose the general It\^o-Taylor schemes for solving MSDEJs.
The authors studied the mean-field It\^o formula
 and proposed the general It\^o-Taylor schemes for MSDEs in \cite{SYZ2017}.
By using the mean-field It\^o formula,
we first develop the It\^o formula for MSDEJs, then based on which, we construct the
It\^o-Taylor expansion for MSDEJs and further propose the It\^o-Taylor schemes of
strong order $\gamma$ and weak order $\eta$ for solving MSDEJs.
Taking $\gamma = 0.5,1.0$ and $\eta = 2.0$,
we obtain the Euler scheme, the strong order $1.0$ Taylor scheme,
and the weak order $2.0$ Taylor scheme, respectively.
Moreover, the rigorous error estimates
indicate that the order of strong convergence
of the strong order $\gamma$ Taylor scheme
is $\gamma$ and the order of weak convergence of the weak order $\eta$ Taylor scheme is $\eta$.
Some numerical tests are carried out
to show the efficiency and the accuracy of the proposed schemes for solving MSDEJs
and to verify our theoretical conclusions.
The numerical results are consistent with our theoretical ones and
show that the efficiency of the proposed schemes
depends on the level of the intensity of the Poisson random measure.

\section{Preliminaries}\label{sec:preliminary}

\subsection{Existence and uniqueness of solution of MSDEJs}\label{subsec:solution}
In this subsection, we state a standard result on the
existence and uniqueness of the strong solutions of the MSDEJ \eqref{MSDEJ}.
For this end, we set the following assumptions on $b$, $\sigma$  and $c$.
\vspace{0.02cm}
\begin{enumerate}
\item[${\rm(A1)}$]
$b(\cdot,x',x)$, $\sigma(\cdot,x',x)$
and $c(\cdot,x',x,e)$ are deterministic continuous processes,
for any fixed $(x',x,e)\in\R^d\times\R^d\times\mathsf{E}$.
\item[${\rm(A2)}$]
There exists a positive constant $L$ such that
\begin{align*}
&|b(t,x',x)-b(t,y',y)|+|\sigma(t,x',x)-\sigma(t,y',y)|\\
\leq& L(|x'-y'|+|x-y|),
\quad{\rm~~ for~ all \;\;t\in[0, T]\;\; and\;\; x,x',y,y'\in \R^d.}
\end{align*}	
\item[${\rm(A3)}$]
There exists a function
$\rho:\mathsf{E}\rightarrow \R^{+}$ satisfying $\int_{\mathsf{E}}\rho^2(e)\lambda(de)< +\infty$,
such that
\begin{align*}
&|c(t,x',x,e) - c(t,y',y,e)|\leq  \rho(e)(|x'-y'| + |x-y|),\\
&|c(t,0,0,e)|\le \rho(e),\
\quad{\rm~~ for~ all \;\;t\in[0, T],\;\; x,x',y,y'\in \R^d\;\;\text{and}\;\;e\in\msE.}
\end{align*}
\item[${\rm(A4)}$]
There exists a constant $K>0$
such that
\begin{align*}
&|b(t,x',x)| + |\sigma(t,x',x)| \leq K(1+|x|+|x'|),\\
&|c(t,x',x,e)|\le \rho(e)(1+|x'|+|x|),
\end{align*}
for all $t\in[0, T]$, $x,x'\in \R^d$, and $e\in\msE$.
\end{enumerate}
\vspace{0.05cm}
%
%

Now we state the existence and uniqueness of the solution of the MSDEJ \eqref{MSDEJ} and some
useful properties in the following theorem \cite{HL2016,Li2017}.

\begin{thm}\label{th0}
Suppose that the coefficients $b$, $\sigma$ and $c$ of the MSDEJ \eqref{MSDEJ}
satisfy the assumptions ${\rm (A1)} - {\rm (A4)}$,
and the initial data $\xi$ and $\xi'$ satisfy
\[\E\big[|\xi|^2+|\xi'|^2\big]<+\infty,\]
then the MSDEJ \eqref{MSDEJ} admits a unique strong solution $X_t^{t_0,\xi}$ on $[t_0,T]$
with
\begin{equation}\label{estimate}
\sup_{t_0\le t\le T}\E\left[\big|X_t^{t_0,\xi}\big|^2\right]<+\infty.
\end{equation}

In addition, for any $p\ge 2$, there exists a $C_p\in \R^+$ such that
for any initial time $t_0\in [0,T]$ and $\mathcal F_{t_0}$
measurable  $\xi_1,\xi_2\in L^p$,
\begin{subequations}\label{con}
\begin{align}
\E\bigg[\sup\limits_{t_0\leq s\leq T}\big|X_s^{t_0,\xi_1}\big|^p\Big|\mathcal{F}_{t_0}\bigg]
&\leq C_p\big(1+|\xi_1|^p\big), \\
\E\bigg[\sup\limits_{t_0\leq s\leq T}\big|X_s^{t_0,\xi_1}-X_s^{t_0,\xi_2}\big|^p\Big|\mathcal{F}_{t_0}\bigg]
&\leq C_p(|\xi_1-\xi_2|^p), \\
\E\bigg[\sup \limits_{t_0\leq s\leq t_0+\delta}\big|X_s^{t_0,\xi_1}-\xi_1\big|^p\Big|\mathcal{F}_{t_0}\bigg]
&\leq C_p\big(1+|\xi_1|^p\big)\delta,
\end{align}
\end{subequations}
a.s. for all $\delta >0$ with $t_0+\delta\leq T$,
where $X_s^{t_0,\xi_1}$ and $X_s^{t_0,\xi_2}$
are the solutions of \eqref{MSDEJ} with initial conditions
$\xi_1$ and $\xi_2$, respectively.
Here the constant $C_p$ in \eqref{con} only depends on $L$, $K$ and $\rho(e)$.
\end{thm}

\subsection{The Markov property}

In this subsection, we present the Markov property of
the solutions of MSDEJs, which will play a key role in
our error estimates.
%
For simplicity, we let $t_0=0$ and denote by
$X_0=\xi$ and $X_t=X_t^{0,\xi}$, then the MSDEJ \eqref{MSDEJ} becomes
\begin{equation}\label{MSDEJsimp}
\begin{aligned}
X_t= X_0 &+\int_0^t\E\big[b(s,X_s^{0,\xi'},x\big]\big|_{x=X_s} ds
+ \int_0^t\E\big[\sigma(s,X_s^{0,\xi'},x)\big]\big|_{x=X_s} dW_s \\
& + \int_0^t\int_{\msE} \E\big[c(s,X_{s-}^{0,\xi'},x,e)\big]\big|_{x=X_{s-}}\mu(de,ds).
\end{aligned}
\end{equation}
Let $X_t^{s,x}$ be the solution of the MSDEJ \eqref{MSDEJsimp}
starting from the point $(s,x)$, i.e.,
\begin{equation*}
\begin{aligned}
X_t^{s,x}= x &+\int_s^t\E\big[b(r,X_r^{0,\xi'},x)\big]\big|_{x=X_r^{s,x}} dr
+ \int_s^t\E\big[\sigma(r,X_r^{0,\xi'},x)\big]\big|_{x=X_r^{s,x}} dW_r \\
& + \int_s^t\int_{\msE} \E\big[c(r,X_{r-}^{0,\xi'},x,e)\big]\big|_{x=X_{r-}^{s,x}}\mu(de,dr)
\end{aligned}
\end{equation*}
for $0\le s\le t\le T$. Then by Theorem \ref{th0}, we have
\begin{equation}\label{eq:eqsolu}
\begin{aligned}
X_t^{s,X_s}=X_t.
\end{aligned}
\end{equation}

Now we state the Markov property of the solutions of MSDEJs as below.
\begin{thm}[The Markov property]\label{th:markov}
Let $f:\R^d\rightarrow\R$ be a bounded Borel measurable function.
Then for the solution $X_t$ of \eqref{MSDEJsimp}, it holds that
\begin{equation}\label{Markov}
\E\left[f(X_t)\big|\F_s\right]
=\E\left[f\big(X_t^{s,X_s}\big)\right]
=\E\left[f\big(X_t^{s,y}\big)\right]\big|_{y=X_s},
\end{equation}
where $0\le s\le t\le T$.
\end{thm}
\begin{proof}
By using the relationship \eqref{eq:eqsolu},
the proof of Theorem \ref{th:markov} is similar
to that of Theorem 7.1.2 in \cite{o}.
So we omit it here.
\end{proof}
\subsection{The equivalent form of MSDEJs}\label{subsec:martingale}
By using the relationship $\tilde \mu(de,dt)=\mu(de,dt)-\lambda(de)dt$,
the MSDEJ \eqref{MSDEJ} can be written as
\begin{equation}
\begin{aligned}\label{MSDEJ-martingale}
X_t^{t_0,\xi} =\,& \xi + \int_{t_0}^t\E\big[\tilde{b}(s,X_s^{t_0,\xi'},x)\big]\big|_{x=X_s^{t_0,\xi}} ds
+ \int_{t_0}^t\E\big[\sigma(s,X_s^{t_0,\xi'},x)\big]\big|_{x=X_s^{t_0,\xi}}  dW_s \\
& \;\;+ \int_{t_0}^t\int_{\msE} \E\big[c(s,X_{s-}^{t_0,\xi'},x,e)\big]\big|_{x=X_{s-}^{t_0,\xi}}
\tilde{\mu}(de,ds),
\end{aligned}
\end{equation}
where the compensated drift coefficient $\tilde{b}$ is defined by
\begin{equation}\label{tildedrift}
\tilde{b}(t,x',x)=b(t,x',x)+\int_{\msE}c(t,x',x,e)\lambda(de).
\end{equation}
Note that by \eqref{tildedrift} and the
assumptions ${\rm (A1)} - {\rm (A4)}$, we can conclude that $\tilde{b}$ satisfies the Lipschitz condition
\[|\tilde{b}(t,x',x)-\tilde{b}(t,y',y)|\le C(|x'-y'| + |x-y|),\]
as well as the linear growth condition
\[|\tilde{b}(t,x',x)|\le C(1+|x'|+|x|),\]
for $t\in[0,T]$ and $x,y,x',y'\in\R^d$.
Here $C$ is a constant depending on $L$, $K$ and $\rho(e)$.

Based on the two equivalent forms of the MSDEJs \eqref{MSDEJ} and \eqref{MSDEJ-martingale},
we will derive two different types of It\^o-Taylor schemes for solving MSDEJs.
\section{The It\^o formula and It\^o-Taylor expansion}\label{sec:ito}

In this section, we develop the It\^o formula
and It\^o-Taylor expansion for MSDEJs,
which are the foundation for proposing the
It\^o-Taylor schemes for MSDEJs.

\subsection{It\^o's formula for MSDEJs}
In this subsection, based on the mean-field It\^o formula \cite{SYZ2017},
we rigorously prove the It\^o's formula for MSDEJs.

Let $X_t$ be a $d$-dimensional It\^o process
satisfying the MSDE
\begin{equation}\label{itxx}
dX_t= b^{\beta}(t,X_t)dt+\sigma^{\beta}(t,X_t) dW_t, \quad 0\le t\le T
\end{equation}
with $\beta_t$ a $d$-dimensional It\^o process defined as
\begin{equation}\label{itppp}
d\beta_t=\psi_tdt+\varphi_tdW_t,
\end{equation}
where $\psi_t$ and $\varphi_t$
are two progressively measurable processes such that
$\int_{0}^{T}|\psi_t|dt<+\infty$ and $\int_{0}^{T}{\rm Tr}[\varphi_s\varphi_s^{\top}]dt<+\infty$.
Here $\text{Tr}[A]$ denotes the trace of a matrix $A$.

Moreover, we define $f^{\beta}(t,x)$ and $g^{\beta}(t,x,e)$ by
\begin{equation*}\label{itf}
f^{\beta}(t,x)=\E\left[f(t,\beta_t,x)\right],\;\;\;
g^{\beta}(t,x,e)=\E\left[g(t,\beta_t,x,e)\right],
\end{equation*}
for functions $f(t,x',x):\R^+\times \R^d\times \R^d\rightarrow \R$
and $g(t,x,x',e):\R^+\times \R^d\times \R^d\times\msE\rightarrow \R$.
Then we have the following It\^o's formula for the MSDE \eqref{itxx}.
\begin{thm}[Mean-field It\^o formula \cite{SYZ2017}]\label{thm2}
Let $X_t$ and $\beta_t$ be $d$-dimensional It\^o processes satisfying \eqref{itxx}
and \eqref{itppp}, respectively,
and function $f=f(t,x',x)\in C^{1,2,2}$.
Then $f^{\beta}(t,X_t)$ is an It\^o process and satisfies
\begin{equation}\label{mit}
\begin{aligned}
f^{\beta}(t,X_t)
= f^{\beta}(0,X_0)+\int_{0}^{t}L^0f^{\beta}(s,X_s)ds
+ \int_{0}^{t}\overrightarrow{L}^1 f^{\beta}(s,X_s)dW_s,
\end{aligned}
\end{equation}
where  $L^0$ and $\overrightarrow{L}^1$
are defined by
\begin{equation}\label{def:l0}
\begin{aligned}
L^0f^\beta(s,x)= \;& \frac{\partial f^{\beta}}{\partial s}(s,x)
+\nabla_{x}f^\beta(s,x)b^{\beta}(s,x)\\
%
&+ \frac{1}{2}{\rm Tr}\big[\nabla_{xx}f^\beta(s,x)\big(\sigma^{\beta}(s,x)\big)\big(\sigma^{\beta}(s,x)\big)^\top\big],\\
\overrightarrow{L}^1f^\beta(s,x) = \;& \nabla_{x}f^\beta(s,x)\sigma^{\beta}(s,x)
=\big(L^1 f^\beta(s,x), \ldots, L^m f^\beta(s,x)\big) \\
\end{aligned}
\end{equation}
with
\begin{align*}
&
L^j f^\beta(t,x)=\sum\limits_{k=1}^{d}
\frac{\partial f^\beta}{\partial x^k}(t,x)\sigma_{kj}^{\beta}(t,x),\quad j=1,\dots,m, \notag \\
\label{ft}
& \frac{\partial f^\beta}{\partial s}(s,x)
=\E\bigg[\frac{\partial f}{\partial s}(s,\beta_s,x)+\nabla_{x'}f(s,\beta_s,x)\psi_s
+\frac{1}{2}{\rm Tr}\big[f_{x'x'}(s,\beta_s,x)\varphi_s\varphi_s^\top\big]\bigg],
\end{align*}
and
$$\begin{aligned}	
\nabla_{x} f^\beta(s,x) = \E\big[\nabla_{x}f(s,\beta_s,x)\big],\;\;\;\;\;\;
f_{xx}^\beta(s,x) = \E\big[f_{xx} (s,\beta_s,x)\big],\qquad\qquad
\end{aligned}$$
where
$\nabla_xf=\Big(\frac{\partial f}{\partial x_1},\ldots,\frac{\partial f}{\partial x_d}\Big)$
is a $d$-dimensional row vector,
$f_{xx} = \Big(\frac{\partial^2f}{\partial x_i\partial x_j}\Big)_{d\times d}$
a $d\times d$ matrix.	
\end{thm}

Now let $X_t$ be a $d$-dimensional It\^o process with jumps
satisfying the MSDEJ
\begin{equation}\label{ito-msdej}
dX_t= b^{\beta}(t,X_t)dt+\sigma^{\beta}(t,X_t) dW_t
 + \int_{\msE}c^{\beta}(t-,X_{t-},e)\mu(de,dt), \quad 0\le t\le T
\end{equation}
with $\beta_t$ a $d$-dimensional It\^o process with jumps defined by
\begin{equation}\label{itpppj}
d\beta_t=\psi_tdt+\varphi_tdW_t+\int_{\msE}h_t\mu(de,dt),
\end{equation}
where $h_t$ is a progressively measurable process such that
$\int_{\msE}|h_t|\lambda(de)<+\infty$.
Here by the definition of $g^{\beta}(t,x,e)$, we have
\begin{equation}\label{definitionscc}
\left\{\begin{aligned}
&c^{\beta}(t,x,e)\;\;\,=\E\left[c(t,\beta_{t},x,e)\right],\\
&c^{\beta}(t-,x,e)=\E\left[c(t,\beta_{t-},x,e)\right].
\end{aligned}\right.
\end{equation}
\begin{rem}\label{remark1}
By the property of the Poisson process $N_t=\mu(\msE\times[0,t])$,
for any fixed $t\in[0,T]$, with probability $1$,
$t$ is not a jump time \cite{CT2004}.
Then by \eqref{itpppj}, it holds that
$\beta_t=\beta_{t-}$ a.s., which leads to
$$\E\left[c(t,\beta_{t},x,e)\right]
=\E\left[c(t,\beta_{t-},x,e)\right],$$
that is,
\begin{equation}\label{eq:mean}
c^{\beta}(t,x,e)=c^{\beta}(t-,x,e).
\end{equation}
\end{rem}

Now  we state the It\^o's formula for MSDEJs in the following theorem.

\begin{thm}[Mean-field It\^o formula with jumps]\label{thm3}
Let $X_t$ and $\beta_t$ be two $d$-dimensional It\^o processes
with jumps defined by \eqref{ito-msdej}
and \eqref{itpppj}, respectively,
and function $f=f(t,x',x)\in C^{1,2,2}$.
Then $f^{\beta}(t,X_t)$ is an It\^o process with jumps and satisfies
\begin{equation}\label{mitj}
\begin{aligned}
f^{\beta}(t,X_t)
= f^{\beta}(0,X_0)&\;+\int_{0}^{t}L^0f^{\beta}(s,X_s)ds
+ \int_{0}^{t}\overrightarrow{L}^1 f^{\beta}(s,X_s)dW_s\\
&\;+\int_0^t\int_{\msE}L_e^{-1}f^{\beta}(s,X_{s-})\mu(de,ds),
\end{aligned}
\end{equation}
where  $L^0$ and $\overrightarrow{L}^1$
are defined by \eqref{def:l0}, and
\begin{equation}\label{def:11}
L_e^{-1}f^{\beta}(s,x) = f^{\beta}(s,x+c^{\beta}(s,x,e))-
f^{\beta}(s,x).
\end{equation}
\end{thm}
\begin{proof}
For simplicity, we consider the case $d=m=q=1$.
The general case can be obtained similarly.

Assume that the Poisson random measure $\mu$
generates a sequence of pairs $\{(\tau_i,Y_i),\\i=1,2,\dots,N_t\}$,
where $N_t=\mu(\msE\times[0,t])$
represents the
total number of jumps of $\mu$ up to time $t$,
and ($\tau_i,Y_i$) are the $i$th jump time and jump size, respectively.
Then we can write the MSDEJ \eqref{ito-msdej} as
\begin{equation}\label{ito-msdej-sum}
X_t=X_0+\int_0^t b^{\beta}(s,X_s)ds+\int_{0}^t\sigma^{\beta}(s,X_s) dW_s
 + \sum_{i=1}^{N_t}c^{\beta}(\tau_i,X_{\tau_i-},Y_i),
\end{equation}
where $c^{\beta}(\tau_i,X_{\tau_i-},Y_i)=\E\big[c(\tau_i,\beta_{\tau_i},x,e)\big]\big|_{(x,e)=(X_{\tau_i-},Y_i)}$.

Let $\tau_0=0$ and $\tau_{N_t+1}=t$, and we have
\begin{equation}\label{ito-msdejs}
\begin{aligned}
f^{\beta}(t,X_t)-f^{\beta}(0,X_0)=&\;\sum_{i=0}^{N_t}\left(f^{\beta}(\tau_{i+1},X_{\tau_{i+1}}) -f^{\beta}(\tau_i,X_{\tau_{i}})\right)\\
 =&\;\sum_{i=0}^{N_t}\big(f^{\beta}(\tau_{i+1},X_{\tau_{i+1}}) -f^{\beta}(\tau_{i+1}-,X_{\tau_{i+1}-})\big)\\
&+ \sum_{i=0}^{N_t}\big(f^{\beta}(\tau_{i+1}-,X_{\tau_{i+1}-}) -f^{\beta}(\tau_i,X_{\tau_{i}})\big).
\end{aligned}
\end{equation}
Note that, $X_t$ is a MSDE on each time interval $[\tau_i,\tau_{i+1})$ for $i=0,\dots,N_t$,
then by the mean-field It\^o formula \eqref{mit}, we obtain
\begin{equation}\label{ito-smooth}
\begin{aligned}
f^{\beta}(\tau_{i+1}-,X_{\tau_{i+1}-})-f^{\beta}(\tau_{i},X_{\tau_i})= \int_{\tau_i}^{\tau_{i+1}-}\!\!L^0f^{\beta}(s,X_s)ds
+\int_{\tau_i}^{\tau_{i+1}-}\!\!L^1f^{\beta}(s,X_s)dW_s.
\end{aligned}
\end{equation}
According to Remark \ref{remark1}, at each jump time $\tau_i$, $i=1,\dots,N_t$, $f^{\beta}(t,X_t)$ has a jump 
\begin{equation}\label{ito-jumppart}
\begin{aligned}
f^{\beta}(\tau_i,X_{\tau_{i}})-f^{\beta}(\tau_{i}-,X_{\tau_{i}-})=f^{\beta}\big(\tau_i,X_{\tau_{i}-}+c^{\beta}(\tau_i,X_{\tau_i-},Y_i)\big)
-f^{\beta}(\tau_i,X_{\tau_{i}-}).
\end{aligned}
\end{equation}
%
Then by \eqref{ito-msdejs} - \eqref{ito-jumppart} we get
\begin{equation*}\label{ito-msdejs-smooth}
\begin{aligned}
f^{\beta}(t,X_t)-f^{\beta}(0,X_0)
=&\; \sum_{i=0}^{N_t}\big(f^{\beta}(\tau_{i+1}-,X_{\tau_{i+1}-}) -f^{\beta}(\tau_i,X_{\tau_{i}})\big)\\
&\;+\sum_{i=1}^{N_t}\big(f^{\beta}(\tau_{i},X_{\tau_{i}}) -f^{\beta}(\tau_{i}-,X_{\tau_{i}-})\big)\\
=&\;\sum_{i=0}^{N_t}\big(\int_{\tau_i}^{\tau_{i+1}-}L^0f^{\beta}(s,X_s)ds
+\int_{\tau_i}^{\tau_{i+1}-}L^1f^{\beta}(s,X_s)dW_s\big)\\
&\;+\sum_{i=1}^{N_t}\big(f^{\beta}(\tau_i,X_{\tau_{i}-}+c^{\beta}(\tau_i,X_{\tau_i-},Y_i))
-f^{\beta}(\tau_i,X_{\tau_{i}-})\big)\\
=&\;\int_{0}^{t}L^0f^{\beta}(s,X_s)ds
+ \int_{0}^{t}\overrightarrow{L}^1 f^{\beta}(s,X_s)dW_s\\
&\;+\int_0^t\int_{\msE}L_e^{-1}f^{\beta}(s,X_{s-})\mu(de,ds),
\end{aligned}
\end{equation*}
where $L_e^{-1}$ is defined by \eqref{def:11}.
We complete the proof.
\end{proof}

By using the relationship
$\tilde{\mu}(de,dt) = \mu(de,dt)-\lambda(de)dt$,
we  write \eqref{ito-msdej} as
\begin{equation}\label{ito-msdej-com}
dX_t= \tilde{b}^{\beta}(t,X_t)dt+\sigma^{\beta}(t,X_t) dW_t
 + \int_{\msE}c^{\beta}(t,X_{t-},e)\tilde{\mu}(de,dt), \quad 0\le t\le T,
\end{equation}
where $\tilde{b}^{\beta}(t,x)=b^{\beta}(t,x)+\int_{\msE}c^{\beta}(t,x,e)\lambda(de)$.
Based on Theorem \ref{thm3},
we have  the It\^o's formula for the equivalent  MSDEJ \eqref{ito-msdej-com} as below.
\begin{pro}
Let $X_t$ and $\beta_t$ be two $d$-dimensional It\^o processes with jumps defined by
\eqref{ito-msdej-com} and \eqref{itpppj}, respectively,
and function $f=f(t,x',x)\in C^{1,2,2}$.
Then $f^{\beta}(t,X_t)$ is an It\^o process with jumps and satisfies
 \begin{equation}\label{mitj-martingale}
\begin{aligned}
f^{\beta}(t,X_t)
= f^{\beta}(0,X_0)&\;+\int_{0}^{t}\tilde{L}^0f^{\beta}(s,X_s)ds
+ \int_{0}^{t}\overrightarrow{L}^1 f^{\beta}(s,X_s)dW_s\\
&\;+\int_0^t\int_{\msE}L_e^{-1}f^{\beta}(s,X_{s-})\tilde{\mu}(de,ds),
\end{aligned}
\end{equation}
where 
\begin{equation}\label{tildel0}
\tilde{L}^0 f^{\beta}(s,x)= L^0f^{\beta}(s,x) + \int_{\msE}L_e^{-1}f^{\beta}(s,x)\lambda(de).
\end{equation}
\end{pro}

Note that \eqref{mitj-martingale} is an equivalent form
of the It\^o's formula \eqref{mitj} for MSDEJs.
Moreover, when $f$ is independent of $x'$, the It\^o's formulas \eqref{mitj} and \eqref{mitj-martingale}
for MSDEJs reduce to the ones for standard SDEJs \cite{CT2004,PB2010}.
Hence the It\^o's formulas for MSDEJs
can be seen as a generalization of
the ones for SDEJs.
%
%

\subsection{It\^o-Taylor expansion for MSDEJs}\label{c6s23}

In this subsection, by utilizing It\^o's formula,
we construct the It\^o-Taylor expansion for MSDEJs.
To proceed, we introduce
multiple It\^o integrals and
coefficient functions as below.
\subsubsection{Multiple It\^o integrals}
In this subsection, we introduce two types of multiple stochastic integrals.

\begin{enumerate}[(A)]
\vspace{0.15cm}
\item Multi-indices\\
Let $\alpha=(j_1,\cdots,j_l)$ be a multi-index with
$j_i\in\{-1,0,1, \cdots,m\},\; i=1,\cdots, l$.
Set $l(\alpha)=l$ to be the length of $\alpha$,
and let ${\mathcal{M}}$ be the set of all multi-indices, i.e.,
$$\mathcal{M}=\Big\{(j_1,j_2,\cdots,j_l): j_i\in\{-1,0,1,\cdots,m\},~i\in\{1,2,\cdots,l\},\;
l\in\mathcal{N}^+\Big\}\cup\{v\},$$
where $v$ is the multi-index of length zero,
i.e.,\; $l(v) = 0$.
For a given $\alpha\in \mathcal{M}$ with $l(\alpha)\ge1$,
$-\alpha$ and $\alpha-$ are two multi-indices
obtained by deleting the first and the last component of $\alpha$, respectively.
We also denote by
\begin{align*}
n(\alpha)&: ~\text{the number of the components of $\alpha$ equal to 0},\\
s(\alpha)\,&: ~\text{the number of the components of $\alpha$ equal to -1}.
\end{align*}

Moreover, for a given $\alpha\in\mathcal{M}$, let
$e=(e_1,\dots,e_{s(\alpha)})$ denote a vector $e\in\msE^{s(\alpha)}$.
\vspace{0.1cm}
\item Multiple integrals\\
For a given $\alpha\in\mathcal{M}$, we define the multiple It\^o integral operator $I_\alpha$ on
the adapted right continuous processes $\{f=f(t,e_1,\dots,e_{s(\alpha)}),t\geq 0\}$
with left  limits by
\begin{equation}\label{iti}
I_\alpha[f(\cdot)]_{\rho,\tau}:=
\begin{cases}
f(\tau), & ~l=0,\\
\int_{\rho}^{\tau}I_{\alpha-}[f(\cdot)]_{\rho,s}ds, & ~l\ge 1 ~\text{and}~j_l=0,\\
\int_{\rho}^{\tau}I_{\alpha-}[f(\cdot)]_{\rho,s}dW_s^{j_l}, & ~l\ge 1 ~\text{and}~j_l\ge 1,\\
\int_{\rho}^{\tau}\!\int_{\msE}I_{\alpha-}[f(\cdot)]_{\rho,s-}\mu(de_{s(\alpha)},ds), & ~l\ge 1 ~\text{and}~j_l= -1,
\end{cases}
\end{equation}
where $\rho$ and $\tau$ are two stopping times satisfying
$0\leq \rho\leq \tau \leq T$, a.s.
and all the integrals exist.
For instance,
\begin{align*}
&I_v[f(\cdot)]_{0,t}=f(t), \quad I_{(-1)}[f(\cdot)]_{0,t}=\int_{0}^{t}\!\int_{\msE}f(s-,e)\mu(de,ds),\\
&I_{(0)}[f(\cdot)]_{0,t}=\int_{0}^{t}f(s)ds,\quad I_{(1)}[f(\cdot)]_{0,t}=\int_{0}^{t}f(s)dW_s^1,\\
&I_{(-1,-1)}[f(\cdot)]_{0,t}=\int_{0}^{t}\!\int_{\msE}\int_0^{s_2-}\!\!\!\int_{\msE}f(s_1-,e_1,e_2)\mu(de_1,ds_1)\mu(de_2,ds_2).
\end{align*}
%
\item Compensated multiple integrals\\
Replace $\mu$ with $\tilde{\mu}$ in \eqref{iti},
and we get the compensated multiple It\^o integral
\begin{equation*}\label{iticom}
\tilde{I}_\alpha[f(\cdot)]_{\rho,\tau}:=
\begin{cases}
f(\tau), & ~l=0,\\
\int_{\rho}^{\tau}\tilde{I}_{\alpha-}[f(\cdot)]_{\rho,s}ds, & ~l\ge 1 ~\text{and}~j_l=0,\\
\int_{\rho}^{\tau}\tilde{I}_{\alpha-}[f(\cdot)]_{\rho,s}dW_s^{j_l}, & ~l\ge 1 ~\text{and}~j_l\ge 1,\\
\int_{\rho}^{\tau}\!\int_{\msE}\tilde{I}_{\alpha-}[f(\cdot)]_{\rho,s-}\tilde{\mu}(de_{s(\alpha)},ds), & ~l\ge 1 ~\text{and}~j_l= -1,
\end{cases}
\end{equation*}
where $f(\cdot) = f(\cdot,e_1,\dots,e_{s(\alpha)})$.
For instance,
\begin{align*}
\tilde{I}_{\alpha}[f(\cdot)]_{\rho,\tau}&=I_{\alpha}[f(\cdot)]_{\rho,\tau}, \;\;\;\;\;\;\;\;\;\text{when}\;\;\;s(\alpha)=0,\\
\tilde{I}_{(1,-1)}[f(\cdot)]_{0,t}\;\,&=\int_{0}^{t}\!\int_{\msE}\int_{0}^{s_2-}\!f(s_1,e)dW_{s_1}^1\tilde{\mu}(de,ds_2).
\end{align*}
\end{enumerate}
\subsubsection{Coefficient functions}
For a given function
$$f(t,x',x,e):\R^+\times \R^d\times\R^d\times\msE\rightarrow\R,$$
by \eqref{def:l0}
and \eqref{def:11}, we have
 \begin{equation}\label{def:l2}
\begin{aligned}
L^0f^\beta(s,x,e)= \;& \frac{\partial f^{\beta}}{\partial s}(s,x,e)
+\sum\limits_{k=1}^{d}b_k^{\beta}(s,x)
\frac{\partial f^\beta}{\partial x^k}(t,x,e)\\
%
&+ \frac{1}{2}\sum\limits_{i,k=1}^{d}\sum_{j=1}^d
\sigma_{ij}^{\beta}(t,x)\sigma_{kj}^{\beta}(s,x)
\frac{\partial^2 f^\beta}{\partial x^i\partial x^k}(t,x,e),
\end{aligned}
\end{equation}
with
\begin{align*}	
\frac{\partial f^\beta}{\partial s}(s,x,e)
=\E\bigg[\frac{\partial f}{\partial s}(s,\beta_s,x,e)+\nabla_{x'}f(s,\beta_s,x,e)\psi_s
+\frac{1}{2}{\rm Tr}\big[f_{x'x'}(s,\beta_s,x,e)\varphi_s\varphi_s^\top\big]\bigg],
\end{align*}
and
\begin{equation}\label{def:l3}
\begin{aligned}
&L^j f^\beta(t,x,e)=\sum\limits_{k=1}^{d}
\frac{\partial f^\beta}{\partial x^k}(t,x,e)\sigma_{kj}^{\beta}(t,x),\quad j=1,\dots,m,\\
&L_{e_2}^{-1}f^{\beta}(s,x,e_1)=f^{\beta}(s,x+c^{\beta}(s,x,e_2),e_1)-f^{\beta}(s,x,e_1).
\end{aligned}\end{equation}
Then by \eqref{tildel0}, we obtain
\begin{equation}\label{def:l4}
\tilde{L}^{0}f^{\beta}(s,x,e)=L^0f^{\beta}(s,x,e)+\int_{\msE}L_{e_1}^{-1}f^{\beta}(t,x,e)\lambda(de_1).
\end{equation}

Based on \eqref{def:l2} - \eqref{def:l4}, we present the following two types of coefficient functions
and hierarchical and remainder sets.

\begin{enumerate}[(C)]
\item It\^o coefficient functions\\
For a given $\alpha = (j_1,\cdots,j_l)\in\mathcal{M}$
and a smooth function $f(t,x',x)$,
we  define the coefficient function $f_\alpha^\beta$ by
\begin{equation}\label{falpha}
f_\alpha^\beta(t,x,e):=
\begin{cases}
f^\beta(t,x), & l=0,\\
L^{j_1}f_{-\alpha}^\beta(t,x,e_1,\dots,e_{s(-\alpha)}), &l\ge 1 ~\text{and}~j_1\ge 0,\\
L_{e_{s(\alpha)}}^{-1}f_{-\alpha}^\beta(t,x,e_1,\dots,e_{s(-\alpha)}), &l\ge 1 ~\text{and}~j_1=-1,
\end{cases}
\end{equation}
where $e=(e_1,\dots,e_{s(\alpha)})\in\msE^{s(\alpha)}$.
The dependence on $e$ in \eqref{falpha} is introduced by the repeated application of
the operator $L_e^{-1}$ in \eqref{def:l3}.
%
Take $m=d=q=1$ and let $f(t,x',x)=x$, then we can deduce the following examples
$$\begin{aligned}
&f_{(0)}^\beta(t,x) =b^\beta(t,x),\quad f_{(1)}^\beta(t,x)=\sigma^\beta(t,x),
\quad f_{(-1)}^\beta(t,x,e)=c^{\beta}(t,x,e), \\
%
&f_{(-1,-1)}^\beta(t,x,e_1,e_2) =L_{e_2}^{-1}c^\beta(t,x,e_1)=c^{\beta}(t+,x+c^{\beta}(t,x,e_2),e_1)-c^{\beta}(t,x,e_1).
\end{aligned}$$
\item[(D)] Compensated It\^o coefficient functions\\
%
By replacing $L^0$ with $\tilde{L^0}$ in \eqref{falpha},
 we get the compensated It\^o coefficient functions
\begin{equation}\label{tildefalpha}
\tilde{f}_\alpha^\beta(t,x,e):=
\begin{cases}
f^\beta(t,x), & l=0,\\
\tilde{L}^{0}\tilde{f}_{-\alpha}^\beta(t,x,e_1,\dots,e_{s(-\alpha)}), &l\ge 1 ~\text{and}~j_l= 0,\\
L^{j_l}\tilde{f}_{-\alpha}^\beta(t,x,e_1,\dots,e_{s(-\alpha)}), &l\ge 1 ~\text{and}~j_l\ge 1,\\
L_{e_{s(\alpha)}}^{-1}\tilde{f}_{-\alpha}^\beta(t,x,e_1,\dots,e_{s(-\alpha)}), &l\ge 1 ~\text{and}~j_l=-1,
\end{cases}
\end{equation}

For instance, let $f(t,x',x)=x$ and  we have
$$\begin{aligned}
\tilde{f}_{\alpha}^{\beta}(t,x,e)=&\;f_{\alpha}^{\beta}(t,x,e),\;\;\;\;\;\;\;\;\;\qquad\text{when}\;\;\;n(\alpha)=0,\\
\tilde{f}_{(-1,0)}^\beta(t,x,e) =&\;L_e^{-1}\tilde{b}^\beta(t,x)=b^{\beta}(t,x+c^{\beta}(t,x,e))-b^{\beta}(t,x)\\
&\;\,+\int_{\msE}\big(c^{\beta}(t,x+c^{\beta}(t,x,e),e_1)-c^{\beta}(t,x,e_1)\big)\lambda(de_1).
\end{aligned}$$

Here we have assumed that the functions $b$, $\sigma$, $c$ and $f$ satisfy all the smoothness and
integrability conditions needed in the definitions of \eqref{falpha} and \eqref{tildefalpha}.
\item[(E)] Hierarchical and remainder sets\\
We call a subset $\mathcal{A}\subset \mathcal{M}$
a hierarchical set if it satisfies
$$
\mathcal{A}\neq \varnothing, \;\;\; \sup_{\alpha\in \mathcal{A}}l(\alpha)<\infty,\;\
\text{ and }
-\alpha\in\mathcal{A}\;\; \text{for each}~\alpha\in\mathcal{A}\backslash \{v\};
$$
and its remainder set $\mathcal B(\mathcal A)$ is defined by
$$ \mathcal{B}(\mathcal{A})=  \{\alpha \in\mathcal{M}\backslash \mathcal{A}: -\alpha \in \mathcal{A}\}. $$
Take $m=1$ for instance and we give two hierarchical sets
\begin{align*}
\mathcal{A}_{0}=&\;\{v,(-1),(0),(1)\},\\
\mathcal{A}_{1}=&\;\{v,(-1),(0),(1),(1,1),(1,-1),(-1,1),(-1,-1)\},
\end{align*}
and their remainder sets are
$$\begin{aligned}
\mathcal{B(A}_0)&=  \left\{(-1,-1),(0,-1),(1,-1),(-1,0),(0,0),(1,0),(-1,1),(0,1),(1,1)\right\},\\
\mathcal{B(A}_1)&= \{(0,-1),(-1,0),(0,0),(1,0),(0,1),(-1,1,1),(0,1,1),\\
&\quad\;\;\,\,(-1,1,-1),(0,1,-1),(1,1,-1),(-1,-1,1),(0,-1,1),\\
&\quad\;\;\,\,(-1,1,1),(1,-1,1),(-1,-1,-1),(0,-1,-1),(1,-1,-1)\}.
\end{aligned}$$
\end{enumerate}
\subsubsection{The It\^o-Taylor expansion}
In this subsection, by using the It\^o's formula \eqref{mitj} and \eqref{mitj-martingale} for MSDEJs,
we present the It\^o-Taylor expansions of $$f^\beta(t,X_{t})=\E[f(t,\beta_t,x]\big|_{x=X_t}$$
for the solution $X_t$ of the MSDEJ \eqref{ito-msdej} with
$\beta_t$ defined by
\eqref{itpppj} satisfying
$$\E\left[\int_{0}^{T}\!\int_{\msE}|\beta_t|^2\lambda(de)dt\right] <\infty.$$

Now we state the It\^o-Taylor expansions for MSDEJs in the following theorem.
\begin{thm}\label{th33} Let $\rho$ and $\tau$ be two stopping times with
$0\leq \rho\leq \tau\leq T$, a.s..
Then for a given hierarchical set $\mathcal{A}\subset\mathcal{M}$
and a function $f(t,x',x): \R^+\times \R^d \times \R^d\rightarrow  \R$,
we have the It\^o-Taylor expansion
\begin{equation}\label{ite}
f^\beta(\tau,X_\tau)=\sum_{\alpha\in\mathcal{A}}I_\alpha\left[f_\alpha^\beta(\rho,X_\rho)\right]_{\rho,\tau}
+\sum_{\alpha\in\mathcal{B(A)}}I_\alpha\left[f_\alpha^\beta(\cdot,X_{\cdot})\right]_{\rho,\tau},
\end{equation}
and the compensated It\^o-Taylor expansion
\begin{equation}\label{itecompensated}
f^\beta(\tau,X_\tau)=\sum_{\alpha\in\mathcal{A}}\tilde{I}_\alpha\left[\tilde{f}_\alpha^\beta(\rho,X_\rho)\right]_{\rho,\tau}
+\sum_{\alpha\in\mathcal{B(A)}}\tilde{I}_\alpha\left[\tilde{f}_\alpha^\beta(\cdot,X_{\cdot})\right]_{\rho,\tau},
\end{equation}
provided that all of the coefficient functions $f_{\alpha}^\beta$ and $\tilde{f}_{\alpha}^{\beta}$
are well defined and all of the multiple It\^o integrals exist.
\end{thm}
\begin{proof}
By an iterated application of the It\^o's formulas \eqref{mitj} and \eqref{mitj-martingale},
the proof of the above theorem is analogous to the ones of
the It\^o-Taylor expansions for standard SDEs \cite{kp} and SDEJs \cite{PB2010}.
So we omit it here.
\end{proof}

We list some remarks for the It\^o-Taylor expansions for MSDEJs as below.
\vspace{0.1cm}
\begin{itemize}
\item For notational simplicity, we have suppressed the dependence on $e\in\msE^{s(\alpha)}$
in the coefficients $f_{\alpha}$ and $\tilde{f}_{\alpha}$ in \eqref{ite} and \eqref{itecompensated}.

\item
When $f$ is independent of $x'$, the It\^o-Taylor expansions \eqref{ite} and \eqref{itecompensated} for MSDEJs
reduce to the ones for standard SDEJs \cite{PB2010}.
Hence, the It\^o-Taylor expansions for MSDEJs
can be seen as a generalization of
the ones for SDEJs.
\end{itemize}

\section{It\^o-Taylor schemes for MSDEJs}\label{sec:sch}

In this section, based on the It\^o-Taylor expansions \eqref{ite} and \eqref{itecompensated},
we propose the general It\^o-Taylor schemes for solving the MSDEJ \eqref{MSDEJ}.

Without loss of generality, we let $t_0=0$ and $\xi=\xi'$ in \eqref{MSDEJ}.
Then omit the superscript $^{t_0,\xi}$, and  we get
\begin{equation}
\begin{aligned}\label{MSDEJeq}
X_t =X_{0} &+ \int_{0}^t\E\big[b(s,X_s,x)\big]\big|_{x=X_s} ds
+ \int_{0}^t\E\big[\sigma(s,X_s,x)\big]\big|_{x=X_s}  dW_s \\
& + \int_{0}^t\int_{\msE} \E\big[c(s,X_{s},x,e)\big]\big|_{x=X_{s-}}\mu(de,ds).
\end{aligned}
\end{equation}
Note that the MSDEJ \eqref{MSDEJeq} has an equivalent form
\begin{equation*}
\begin{aligned}
X_t =X_{0} &+ \int_{0}^t\E\big[\tilde{b}(s,X_s,x)\big]\big|_{x=X_s} ds
+ \int_{0}^t\E\big[\sigma(s,X_s,x)\big]\big|_{x=X_s}  dW_s \\
& + \int_{0}^t\int_{\msE} \E\big[c(s,X_{s},x,e)\big]\big|_{x=X_{s-}}\tilde{\mu}(de,ds),
\end{aligned}
\end{equation*}
where $\tilde{b}$ is defined by \eqref{tildedrift}.

By choosing different hierarchical sets $\mathcal{A}$
in the It\^o-Taylor expansion \eqref{ite},
we shall derive the two types of strong order $\gamma$ and weak order $\eta$
It\^o-Taylor schemes for solving the MSDEJ \eqref{MSDEJeq}.
%
%
To this end, we take a uniform time partition on $[0,T]$:
$$0=t_0<t_1<\cdots<t_{N}=T,$$
where $t_{k+1}$ is $\mathcal{F}_{t_k}$-measurable for $k=0,1,\dots,N-1$.

Let $X_k$ be the approximation of the solution $X_{t}$
of \eqref{MSDEJeq} at time $t=t_{k}$,
and denote by
\begin{align*}
f^{X_k}(t_{k},X_{k})=&\;\E[f(t_k,X_k,x)]\big|_{x=X_k},\\
g^{X_k}(t_{k},X_{k},e)=&\;\E[g(t_k,X_k,x,e)]\big|_{x=X_k},
\end{align*}
for $f(t,x',x): \R^+\times \R^d \times \R^d\rightarrow \R$
and $g(t,x',x,e): \R^+\times \R^d \times \R^d\times\msE\rightarrow \R$.

\subsection{Strong It\^o-Taylor schemes}\label{subsec:strong}
To construct the strong It\^o-Taylor schemes for the MSDEJ \eqref{MSDEJeq},
for $\gamma = 0.5, 1.0,  1.5, \cdots$,
we define the hierarchical set $\mathcal{A}_\gamma$ by
\begin{equation*}\label{hs}
\mathcal{A}_\gamma = \left\{\alpha \in\mathcal{M}: l(\alpha) + n(\alpha) \leq 2\gamma~ ~\text{or}~~  l(\alpha) = n(\alpha) = \gamma + \frac{1}{2}\right\}
\end{equation*}
and denote its remainder set by $\mathcal B(\mathcal A_\gamma)$.	
Take $f(t,x',x)=x$ and let $\beta_t = X_t$, then
by Theorem \ref{th33}, for $k=0,1,\dots,N-1$, we have the It\^o-Taylor expansion	
\begin{equation}\label{mite}
X_{t_{k+1}}=\sum_{\alpha\in\mathcal{A_\gamma}}I_\alpha\big[f_\alpha^X(t_k,X_{t_k})\big]_{t_k,t_{k+1}}
+\sum\limits_{\alpha\in\mathcal{B(A_\gamma)}}I_\alpha\big[f_\alpha^X(\cdot,X_{\cdot})\big]_{t_k,t_{k+1}},
\end{equation}
and the compensated It\^o-Taylor expansion	
\begin{equation}\label{mitecom}
X_{t_{k+1}}=\sum_{\alpha\in\mathcal{A_\gamma}}\tilde{I}_\alpha\big[\tilde{f}_\alpha^X(t_k,X_{t_k})\big]_{t_k,t_{k+1}}
+\sum\limits_{\alpha\in\mathcal{B(A_\gamma)}}\tilde{I}_\alpha\big[\tilde{f}_\alpha^X(\cdot,X_{\cdot})\big]_{t_k,t_{k+1}}.
\end{equation}

By removing the remainder term in \eqref{mite},
we propose the following general strong order
$\gamma$ It\^o-Taylor scheme for solving the MSDEJ \eqref{MSDEJeq}.
\begin{sch}[Strong order $\gamma$ It\^o-Taylor scheme]\label{SSG}
	\begin{equation}\label{roi}
	\begin{aligned}
	X_{k+1}=\sum_{\alpha\in\mathcal{A_\gamma}}I_\alpha\big[f_\alpha^{X_k}(t_k,X_k)\big]_{t_k,t_{k+1}}.
	\end{aligned}
	\end{equation}
\end{sch}

Similarly, we propose the compensated strong order
$\gamma$ It\^o-Taylor scheme.
\begin{sch}[Compensated strong order $\gamma$  It\^o-Taylor scheme]\label{SSGcom}
	\begin{equation}\label{roicom}
	\begin{aligned}
	X_{k+1}=\sum_{\alpha\in\mathcal{A_\gamma}}\tilde{I}_\alpha\big[\tilde{f}_\alpha^{X_k}(t_k,X_k)\big]_{t_k,t_{k+1}}.
	\end{aligned}
	\end{equation}
\end{sch}

Based on Schemes  \ref{SSG} and \ref{SSGcom}, by taking $\gamma=0.5$ and $1.0$,
we will give some specific strong Taylor schemes for MSDEJs in the following subsections.
\subsubsection{The Euler scheme}
Take $\gamma=0.5$ in Scheme \ref{SSG}, and
we have
$$\mathcal{A}_{0.5}=\{v,(-1),(0),(1)\}$$
and
$$\mathcal{B}(\mathcal{A}_{0.5})=  \left\{(-1,-1),(0,-1),(1,-1),(-1,0),(0,0),(1,0),(-1,1),(0,1),(1,1)\right\}.$$
Then by the It\^o-Taylor expansion \eqref{mite}, for $k=0,\dots,N-1$, we obtain
\begin{equation}\label{eulref}
\begin{aligned}
X_{t_{k+1}}=X_{t_k}&+b^{X}(t_k,X_{t_k})\int_{t_k}^{t_{k+1}}ds+\sigma^{X}(t_k,X_{t_k})\int_{t_k}^{t_{k+1}}dW_s\\
&+\int_{t_k}^{t_{k+1}}\int_{\msE}c^{X}(t_k,X_{t_k},e)\mu(de,ds)+R_1,
\end{aligned}
\end{equation}
where
\begin{equation*}
\begin{aligned}
R_1=&\;\int_{t_k}^{t_{k+1}}\!\!\!\int_{t_k}^{s}L^0b^{X}(z,X_z)dzds
+\int_{t_k}^{t_{k+1}}\!\!\!\int_{t_k}^{s}L^1b^{X}(z,X_z)dW_zds\\
&~+\int_{t_k}^{t_{k+1}}\!\!\!\int_{t_k}^{s}L^0\sigma^{X}(z,X_z)dzdW_s
+\int_{t_k}^{t_{k+1}}\!\!\!\int_{t_k}^{s}L^1\sigma^{X}(z,X_z)dW_zdW_s\\
&~+\int_{t_k}^{t_{k+1}}\!\!\!\int_{\msE}\int_{t_k}^{s-}\!\!\!L^0c^{X}(z,X_z,e)dz\mu(de,ds)\\
&~+\int_{t_k}^{t_{k+1}}\!\!\!\int_{t_k}^{s}\int_{\msE}L_e^{-1}b^{X}(z,X_{z-})\mu(de,dz)ds\\
&~+\int_{t_k}^{t_{k+1}}\!\!\!\int_{\msE}\int_{t_k}^{s-}\!\!\!L^1c^{X}(z,X_z,e)dW_z\mu(de,ds)\\
&~+\int_{t_k}^{t_{k+1}}\!\!\!\int_{t_k}^{s}\int_{\msE}L_e^{-1}\sigma^{X}(z,X_{z-})\mu(de,dz)dW_s\\
&~+\int_{t_k}^{t_{k+1}}\!\!\!\int_{\msE}\int_{t_k}^{s-}\!\int_{\msE}L_{e_2}^{-1}c^{X}(z,X_{z-},e_1)\mu(de_1,dz)\mu(de_2,ds).
\end{aligned}
\end{equation*}

Remove the remainder term $R_1$ in \eqref{eulref},
and we get the strong order $0.5$ It\^o-Taylor scheme for solving the MSDEJ \eqref{MSDEJeq}
\begin{equation}\label{eul}
\begin{aligned}
X_{k+1}=&X_{k}+b^{X_k}(t_k,X_k)\Delta t_k
+\sigma^{X_k}(t_k,X_k)\Delta W_{k}+\int_{t_k}^{t_{k+1}}\!\!\!\int_{\msE}c^{X_k}(t_k,X_k,e)\mu(de,dt)\\
=&X_{k}+b^{X_k}(t_k,X_k)\Delta t_k
+\sigma^{X_k}(t_k,X_k)\Delta W_{k}
+\sum_{i=N_{t_k}+1}^{N_{t_{k+1}}}c^{X_k}(t_k,X_k,Y_i),
\end{aligned}
\end{equation}	
which is the so-called Euler scheme.
Here $$\Delta t_k=t_{k+1}-t_k~~~\text{and}~~~\Delta W_k=W_{t_{k+1}}-W_{t_k},$$ and
$N_t=\mu(\msE\times[0,t])$ is a Poisson process counting the number of jumps of $\mu$ up to time $t$
and $(\tau_i,Y_i)$ are the $i$th jump time and jump size.

\subsubsection{The strong order $1.0$ It\^o-Taylor scheme}
Take $\gamma=1$ in Scheme \ref{SSG}, and
we have
 $$\mathcal{A}_{1}=\{v,(-1),(0),(1),(1,1),(1,-1),(-1,1),(-1,-1)\}$$
 and
$$\begin{aligned}
\mathcal{B(A}_1)&= \{(0,-1),(-1,0),(0,0),(1,0),(0,1),(-1,1,1),(0,1,1),\\
&\qquad(-1,1,1),(-1,1,-1),(0,1,-1),(1,1,-1),(-1,-1,1),\\
&\qquad(0,-1,1),(1,-1,1),(-1,-1,-1),(0,-1,-1),(1,-1,-1)\}.
\end{aligned}$$
Then by \eqref{roi},  we get the strong order 1.0 It\^o-Taylor scheme
\begin{equation}\label{mil}
\begin{aligned}
X_{k+1}=X_{k}&+b^{X_k}(t_k,X_k)\Delta t_k
+\sigma^{X_k}(t_k,X_k)\Delta W_{k}+\int_{t_k}^{t_{k+1}}\!\!\!\int_{\msE}c^{X_k}(t_k,X_k,e)\mu(de,dt)\\
&+\int_{t_k}^{t_{k+1}}\!\!\!\int_{t_k}^{s}L^1\sigma^{X_k}(t_k,X_{k})dW_zdW_s\\
&+\int_{t_k}^{t_{k+1}}\!\!\!\int_{\msE}\int_{t_k}^{s-}L^1c^{X_k}(t_k,X_{k},e)dW_z\mu(de,ds)\\
&+\int_{t_k}^{t_{k+1}}\!\!\!\int_{t_k}^{s}\int_{\msE}L_e^{-1}\sigma^{X_k}(t_k,X_{k})\mu(de,dz)dW_s\\
&+\int_{t_k}^{t_{k+1}}\!\!\!\int_{\msE}\int_{t_k}^{s-}\!\!\!\int_{\msE}L_{e_1}^{-1}c^{X_k}(t_k,X_{k},e_2)\mu(de_1,dz)\mu(de_2,ds).
\end{aligned}
\end{equation}

Combining with the It\^o's formula \eqref{mitj} for MSDEJs and the properties of
jump times, the scheme \eqref{mil} can be written as
\vspace{-0.15cm}
\begin{align*}
X_{k+1}=&\;X_{k}+b^{X_k}(t_k,X_k)\Delta t_k
+\sigma^{X_k}(t_k,X_k)\Delta W_{k}+\!\sum_{i=N_{t_k}+1}^{N_{t_{k+1}}}\!c^{X_k}(t_k,X_k,Y_i)\\
&\;\,+\frac12L^1\sigma^{X_k}(t_k,X_{k})\big((\Delta W_k)^2-\Delta t_k\big)
+\!\!\sum_{i=N_{t_k}+1}^{N_{t_{k+1}}}\!L^1c^{X_k}(t_k,X_{k},Y_i)\big(W_{\tau_i}-W_{t_k}\big)\\
&\;\,+\!\sum_{i=N_{t_k}+1}^{N_{t_{k+1}}}\!\Big(\sigma^{X_k}\big(t_k,X_{k}+c^{X_k}(t_k,X_k,Y_i)\big)-
\sigma^{X_k}(t_k,X_k)\Big)\big(W_{t_{k+1}}-W_{\tau_i}\big)\\
&\;\,+\!\sum_{i=N_{t_k}+1}^{N_{t_{k+1}}}\sum_{j=N_{t_k}+1}^{N_{\tau_i-}}\!
\Big(c^{X_k}\big(t_k,X_{k}+c^{X_k}(t_k,X_k,Y_j),Y_i\big)-c^{X_k}(t_k,X_k,Y_i)\Big),
\end{align*}
which is readily applicable for scenario simulation.

Based on Scheme \ref{SSGcom}, we can  get
the compensated Euler scheme and the compensated strong order 1.0 It\^o-Taylor scheme
in the same way.
It is also worth noting that the compensated Euler scheme and
the Euler scheme are the same.

\subsection{Weak It\^o-Taylor schemes}\label{subsec:weak}

To construct the weak It\^o-Taylor schemes for MSDEJs, for $\eta = 1.0,  2.0, \cdots$,
we define the hierarchical set $\Gamma_\eta$ by
\begin{equation}\label{hw}
\Gamma_\eta = \left\{\alpha \in\mathcal{M}: l(\alpha) \leq \eta\right\}
\end{equation}
and denote its remainder set by $\mathcal B(\Gamma_\eta)$.	
Take $f(t,x',x)=x$ and let $\beta_t = X_t$, then
by Theorem \ref{th33}, for $k=0,1,\dots,N-1$, we have
\begin{equation}\label{mitew}
X_{t_{k+1}}=\sum_{\alpha\in\Gamma_\eta}I_\alpha\big[f_\alpha^X(t_k,X_{t_k})\big]_{t_k,t_{k+1}}
+\sum\limits_{\alpha\in\mathcal{B}(\Gamma_\eta)}I_\alpha\big[f_\alpha^X(\cdot,X_{\cdot})\big]_{t_k,t_{k+1}},
\end{equation}
and 
\begin{equation}\label{mitecomw}
X_{t_{k+1}}=\sum_{\alpha\in\Gamma_\eta}\tilde{I}_\alpha\big[\tilde{f}_\alpha^X(t_k,X_{t_k})\big]_{t_k,t_{k+1}}
+\sum\limits_{\alpha\in\mathcal{B}(\Gamma_\eta)}\tilde{I}_\alpha\big[\tilde{f}_\alpha^X(\cdot,X_{\cdot})\big]_{t_k,t_{k+1}}.
\end{equation}

Remove the remainder term in \eqref{mitew},
and we propose the following general weak order
$\eta$ It\^o-Taylor scheme for solving the MSDEJ \eqref{MSDEJeq}.
\begin{sch}[Weak order $\eta$ It\^o-Taylor scheme]\label{SWG}
	\begin{equation}\label{wroi}
	\begin{aligned}
	X_{k+1}=\sum_{\alpha\in\Gamma_\eta}I_\alpha\big[f_\alpha^{X_k}(t_k,X_k)\big]_{t_k,t_{k+1}}.
	\end{aligned}
	\end{equation}
\end{sch}

Similarly, we propose the compensated weak order
$\eta$ It\^o-Taylor scheme.
\begin{sch}[Compensated weak order $\eta$  It\^o-Taylor scheme]\label{SWGcom}
	\begin{equation}\label{wroicom}
	\begin{aligned}
	X_{k+1}=\sum_{\alpha\in\Gamma_\eta}\tilde{I}_\alpha\big[\tilde{f}_\alpha^{X_k}(t_k,X_k)\big]_{t_k,t_{k+1}}.
	\end{aligned}
	\end{equation}
\end{sch}

Based on Schemes  \ref{SWG} and \ref{SWGcom}, by taking $\eta=1.0$ and $2.0$,
we will present some specific weak Taylor schemes for MSDEJs.
\subsubsection{The Euler scheme}

Taking $\eta=1.0$ in \eqref{hw} leads to
 $$\Gamma_{1.0}=\{v,(-1),(0),(1)\}$$
 and
$$\mathcal{B}(\Gamma_{1.0})=  \left\{(-1,-1),(0,-1),(1,-1),(-1,0),(0,0),(1,0),(-1,1),(0,1),(1,1)\right\}.$$
Then by Scheme \ref{SWG}, we get the Euler scheme \eqref{eul}
\begin{equation*}
\begin{aligned}
X_{k+1}=&X_{k}+b^{X_k}(t_k,X_k)\Delta t_k
+\sigma^{X_k}(t_k,X_k)\Delta W_{k}
+\sum_{i=N_{t_k}+1}^{N_{t_{k+1}}}c^{X_k}(t_k,X_k,Y_i),
\end{aligned}
\end{equation*}	
which is also the weak order $1.0$
It\^o-Taylor scheme.

\subsubsection{The weak order 2.0 It\^o-Taylor scheme}	
Taking $\eta=2.0$ in \eqref{hw} gives
$$\begin{aligned}
\Gamma_{2.0}=&\{v,(-1),(0),(1),(1,1),(1,-1),(-1,1),(-1,-1),\\
&\qquad\qquad(0,0),(1,0),(0,1),(-1,0),(0,-1)\},\\
\end{aligned}$$
and
$$\begin{aligned}
&\mathcal{B}(\Gamma_{2.0})= \{(-1,1,1),(0,1,1),(-1,1,1),(-1,1,-1),(0,1,-1),(1,1,-1),(-1,-1,1),\\
&\qquad\qquad\quad(0,-1,1),(1,-1,1),(-1,-1,-1),(0,-1,-1),(1,-1,-1),(-1,0,0),\\
&\qquad\qquad\quad(0,0,0),(1,0,0),(-1,1,0),(0,1,0),(1,1,0),(-1,0,1),(0,0,1),(1,0,1),\\
&\qquad\qquad\quad(-1,-1,0),(0,-1,0),(1,-1,0),(-1,0,-1),(0,0,-1),(1,0,-1)\}.
\end{aligned}$$
Then by Scheme \ref{SWG}, we get the weak order 2.0 It\^o-Taylor scheme
\begin{align}\label{weak2.0}
X_{k+1}=&\;X_{k}+b^{X_k}(t_k,X_k)\Delta t_k
+\sigma^{X_k}(t_k,X_k)\Delta W_{k}+\int_{t_k}^{t_{k+1}}\!\!\!\int_{\msE}c^{X_k}(t_k,X_k,e)\mu(de,dt)\nonumber\\
&\;\,+\int_{t_k}^{t_{k+1}}\!\!\!\int_{t_k}^{s}L^1\sigma^{X_k}(t_k,X_{k})dW_zdW_s\nonumber\\
&\;\,+\int_{t_k}^{t_{k+1}}\!\!\!\int_{\msE}\int_{t_k}^{s-}L^1c^{X_k}(t_k,X_{k},e)dW_z\mu(de,ds)\nonumber\\
&\;\,+\int_{t_k}^{t_{k+1}}\!\!\!\int_{t_k}^{s}\int_{\msE}L_e^{-1}\sigma^{X_k}(t_k,X_{k})\mu(de,dz)dW_s\nonumber\\
&\;\,+\int_{t_k}^{t_{k+1}}\!\!\!\int_{\msE}\int_{t_k}^{s-}\!\!\!\int_{\msE}L_{e_1}^{-1}c^{X_k}(t_k,X_{k},e_2)\mu(de_1,dz)\mu(de_2,ds)\nonumber\\
&\;\,+\int_{t_k}^{t_{k+1}}\!\!\!\int_{t_k}^{s}L^0b^{X_k}(t_k,X_{k})dzds
+\int_{t_k}^{t_{k+1}}\!\!\!\int_{t_k}^{s}L^1b^{X_k}(t_k,X_{k})dW_zds\nonumber\\
&\;\,+\int_{t_k}^{t_{k+1}}\!\!\!\int_{t_k}^{s}L^{0}\sigma^{X_k}(t_k,X_{k})dzdW_s
+\int_{t_k}^{t_{k+1}}\!\!\!\int_{\msE}\int_{t_k}^{s-}\!L^{0}c^{X_k}(t_k,X_{k},e)dz\mu(de,ds)\nonumber\\
&\;\,+\int_{t_k}^{t_{k+1}}\!\!\!\int_{t_k}^{s}\int_{\msE}L_e^{-1}b^{X_k}(t_k,X_{k})\mu(de,dz)ds.
\end{align}

Combining with the It\^o's formula \eqref{mitj} for MSDEJs and the properties of
jump times, we can rewrite \eqref{weak2.0} as
\begin{equation*}
\begin{aligned}
X_{k+1}=&\;X_{k}+b^{X_k}(t_k,X_k)\Delta t_k
+\sigma^{X_k}(t_k,X_k)\Delta W_{k}+\!\sum_{i=N_{t_k}+1}^{N_{t_{k+1}}}\!c^{X_k}(t_k,X_k,Y_i)\\
&\;\,+\frac12L^1\sigma^{X_k}(t_k,X_{k})\big((\Delta W_k)^2-\Delta t_k\big)
+\!\sum_{i=N_{t_k}+1}^{N_{t_{k+1}}}\!L^1c^{X_k}(t_k,X_{k},Y_i)\big(W_{\tau_i}-W_{t_k}\big)\\
&\;\,+\!\sum_{i=N_{t_k}+1}^{N_{t_{k+1}}}\!\Big(\sigma^{X_k}\big(t_k,X_{k}+c^{X}(t_k,X_k,Y_i)\big)-
\sigma^{X_k}(t_k,X_k)\Big)\big(W_{t_{k+1}}-W_{\tau_i}\big)\\
&\;\,+\!\sum_{i=N_{t_k}+1}^{N_{t_{k+1}}}\sum_{j=N_{t_k}+1}^{N_{\tau_i-}}\!
\Big(c^{X_k}\big(t_k,X_{k}+c^{X_k}(t_k,X_k,Y_j),Y_i\big)-c^{X_k}(t_k,X_k,Y_i)\Big)\\
&\;\,+\frac{1}{2}L^0b^{X_k}(t_k,X_k)(\Delta t_k)^2+L^1b^{X_k}(t_k,X_k)\Delta Z_{k}\\
&\;\,+L^0\sigma^{X_k}(t_k,X_k)(\Delta W_k\Delta t_k-\Delta Z_{k})
+\!\sum_{i=N_{t_k}+1}^{N_{t_{k+1}}}\!L^0c^{X_k}(t_k,X_{k},Y_i)\big(\tau_i-t_k\big)\\
&\;\,+\!\sum_{i=N_{t_k}+1}^{N_{t_{k+1}}}\!\Big(b^{X_k}\big(t_k,X_{k}+c^{X_k}(t_k,X_k,Y_i)\big)-
b^{X_k}(t_k,X_k)\Big)\big(t_{k+1}-\tau_i\big),
\end{aligned}
\end{equation*}
where $\Delta Z_k$ is a random variable defined by
$$\Delta Z_k=\int_{t_k}^{t_{k+1}}\!\!\!\int_{t_k}^{s}dW_zds
= \Delta W_k\Delta t_k-\int_{t_k}^{t_{k+1}}\!\!\!\int_{t_k}^{s}dzdW_s.$$

Similarly, based on Scheme \ref{SWGcom}, we can obtain
the compensated weak order 2.0 It\^o-Taylor scheme for solving the MSDEJ \eqref{MSDEJeq}.
\vspace{0.1cm}

Now, we illustrate how to use the proposed It\^o-Taylor schemes to solve the
MSDEJ \eqref{MSDEJ} with different initial values $\xi$ and $\xi'$.
Without loss of generality, we set $(t_0,\xi')=(0,x_0)$ and $(t_0,\xi)=(0,X_0)$.
Let $X_n^{x_0}$ and $X_n^{X_0}$ denote the numerical solutions of the It\^o-Taylor schemes
with initial values of $x_0$ and $X_0$, respectively.
Then take the Euler scheme \eqref{eul} for instance, and we can rewrite it as
\begin{equation}\begin{aligned}\label{eulinitial}
X_{k+1}^{X_0}=X_{k}^{X_0}+&\;b^{X_k^{x_0}}\big(t_k,X_k^{X_0}\big)\Delta t_k
+\sigma^{X_k^{x_0}}\big(t_k,X_k^{X_0}\big)\Delta W_{k}\\
+&\;\sum_{i=N_{t_k}+1}^{N_{t_{k+1}}}c^{X_k^{x_0}}\big(t_k,X_k^{X_0},Y_i\big),
\end{aligned}\end{equation}
where $f^{X_k^{x_0}}(t_k,x,e)$ denote $\E\big[f(t_k,X_k^{x_0},x,e)\big]$
for $f=b,\sigma$ and $c$.
Now we can apply the scheme \eqref{eulinitial} to solve the MSDEJ \eqref{MSDEJ} by the following two procedures
\vspace{0.1cm}
%
\begin{itemize}\item[ i)]
Take $X_0=x_0$ and we solve \eqref{MSDEJ} by the scheme \eqref{eulinitial}
to obtain $\{X_n^{x_0}\}_{n=0}^{N}$;
\vspace{0.1cm}
\item[ ii)]
Based on $\{X_n^{x_0}\}_{n=0}^{N}$, we get $\{X_n^{X_0}\}_{n=0}^{N}$ by using the scheme \eqref{eulinitial} again.
\end{itemize}
\vspace{0.1cm}

\begin{rem}
Note that when using the strong order 1.0 Taylor scheme  \eqref{mil} and
the weak order 2.0 Taylor scheme
\eqref{weak2.0} to solve MSDEJs,
we need the knowledge of the exact locations of jump times on the time interval $[0,T]$.
Hence, the efficiency of the schemes
depends on
the intensity of the Poisson measure $\mu$.
And the readers are referred to \cite{CT2004,PB2010}
for details of sampling the jump times of the Poisson measure $\mu$.
\end{rem}

\section{Error estimates for strong Taylor schemes}\label{sec:err}
\setcounter{section}{5}\setcounter{equation}{0}
In this section,
based on the relationship between the local and global convergence rates,
we shall  prove the error estimates of the strong order $\gamma$
It\^o-Taylor Scheme \ref{SSG} and
the compensated strong order $\gamma$ It\^o-Taylor Scheme \ref{SSGcom}.

\subsection{The general error estimate theorem}\label{c6s41}
Let $\{X_{t,X}(s)\}_{t\leq s\leq T}$ be the solution of
the MSDEJ \eqref{MSDEJeq} starting from the point $(t,X)$, that is
\begin{equation}\label{ree}\begin{aligned}
X_{t,X}(t+h)=X&+\int_{t}^{t+h}b^{X_{t,X}}\big(s,X_{t,X}(s)\big)ds
+\int_{t}^{t+h}\sigma^{X_{t,X}}\big(s,X_{t,X}(s)\big)dW_s\\
&+\int_{t}^{t+h}\int_{\msE}c^{X_{t,X}}\big(s,X_{t,X}(s-),e\big)\mu(de,ds),
\end{aligned}\end{equation}
where $h \in [0,T-t]$.
Let $\bar{X}_{t,X}(t+h)$
be the one-step approximation of $X_{t,X}(t+h)$,
and $\bar{X}_{0,X_0}(t_k)$ is the corresponding solution of the one-step scheme
\begin{equation}\label{onestep}
\bar{X}_{0,X_0}(t_k)=\bar{X}_{t_{k-1},\bar{X}_{0,X_0}(t_{k-1})}(t_k),
\end{equation}
with $\bar{X}_{0,X_0}(0)=X_0$.
%
For simplicity, we denote $X_{0,X_0}(t_k)$ by $X(t_k)$
and $\bar{X}_{0,X_0}(t_k)$ by $\bar{X}_k$. Then the one-step scheme \eqref{onestep}
becomes
\begin{equation}\label{onestep1}
\bar{X}_{k}=\bar{X}_{t_{k-1},\bar X_{k-1}}(t_{k}).
\end{equation}

To present the general error estimate theorem for the one-step scheme \eqref{onestep},
we first give the following two lemmas.

\begin{lem}\label{l1}
Let $X_{t,X}(s)$ and $X_{t,Y}(s)$ be the solutions of the MSDEJ \eqref{ree}
with initial conditions $X_{t,X}(t)=X$ and $X_{t,Y}(t)=Y$, respectively.
Let $Z= X_{t,X}(t+h)-X_{t,Y}(t+h)-(X-Y)$,
then under the assumptions ${\rm (A1)} - {\rm (A4)}$, we have
\begin{subequations}
\begin{align}
&\E\big[|Z|^2\big] \leq C \E\big[|X-Y|^2\big] h,\label{e4}\\
&\E\big[|X_{t,X}(t+h)-X_{t,Y}(t+h)|^2\big] \leq (1+C h)\E\big[|X-Y|^2\big], \label{e3}
\end{align}
\end{subequations}
where $C$ is a positive constant depending on $\lambda(\msE)$,
the Lipschitz constant $L$ and the function $\rho(e)$ in assumption ${\rm (A3)}$.
\end{lem}

\begin{proof}
For any $0\leq\theta\leq h$, based on \eqref{ree}, by the It\^o's formula \eqref{mitj}
for MSDEJs and the It\^o's isometry formula, we get
\begin{equation*}
\begin{aligned}
&\;\E\big[|X_{t,X}(t+\theta)-X_{t,Y}(t+\theta)|^2\big]\\
=&\;\E\big[|X-Y|^2\big]
 +\E\Bigg[\int_t^{t+\theta}\big|\sigma^{X_{t,X}}\big(s,X_{t,X}(s)\big)-\sigma^{X_{t,Y}}\big(s,X_{t,Y}(s)\big)\big|^2ds\Bigg]\\
& \;+2\E\Bigg[\int_{t}^{t+\theta}\big(X_{t,X}(s)-X_{t,Y}(s)\big)\Big(b^{X_{t,X}}\big(s,X_{t,X}(s)\big)
-b^{X_{t,Y}}\big(s,X_{t,Y}(s)\big)\Big)ds\Bigg]\\
&\;+\E\Bigg[\int_t^{t+\theta}\!\!\!\int_{\msE}\Big(\big|X_{t,X}(s-)-X_{t,Y}(s-)+
c^{X_{t,X}}\big(s,X_{t,X}(s-),e\big)\\
&\qquad\qquad\qquad\quad-c^{X_{t,Y}}\big(s,X_{t,Y}(s-),e\big)\big|^2-\big|X_{t,X}(s-)-X_{t,Y}(s-)\big|^2\Big)\mu(de,ds)\Bigg]\\
=&\;\E\big[|X-Y|^2\big]
 +\E\Bigg[\int_t^{t+\theta}\big|\sigma^{X_{t,X}}\big(s,X_{t,X}(s)\big)-\sigma^{X_{t,Y}}\big(s,X_{t,Y}(s)\big)\big|^2ds\Bigg]\\
&\; +2\E\Bigg[\int_{t}^{t+\theta}\big(X_{t,X}(s)-X_{t,Y}(s)\big)\Big(b^{X_{t,X}}\big(s,X_{t,X}(s)\big)
-b^{X_{t,Y}}\big(s,X_{t,Y}(s)\big)\Big)ds\Bigg]\\
&\;+\E\Bigg[\int_t^{t+\theta}\!\!\!\int_{\msE}\Big(\big|c^{X_{t,X}}\big(s,X_{t,X}(s),e\big)-c^{X_{t,Y}}\big(s,X_{t,Y}(s),e\big)\big|^2\\
&\;\qquad+2\big|X_{t,X}(s)-X_{t,Y}(s)\big|
\big|c^{X_{t,X}}\big(s,X_{t,X}(s),e\big)-c^{X_{t,Y}}\big(s,X_{t,Y}(s),e\big)\big|
\Big)\lambda(de)ds\Bigg].
\end{aligned}
\end{equation*}
Then under the assumption ${\rm (A2)}$ and ${\rm (A3)}$,
we deduce
\begin{align*}
&\;\E\Big[\big|X_{t,X}(t+\theta)-X_{t,Y}(t+\theta)\big|^2\Big]\\
\leq&\; \E\big[|X-Y|^2\big] +4L^2\int_t^{t+\theta}\E\Big[\big|X_{t,X}(s)-X_{t,Y}(s)\big|^2\Big]ds\\
&\;\qquad\qquad\quad\;\,+4L\int_{t}^{t+\theta}\E\Big[\big|X_{t,X}(s)-X_{t,Y}(s)\big|^2\Big]ds\\	
&\;\qquad\qquad\quad\;\,+8K_1\int_t^{t+\theta}\E\Big[\big|X_{t,X}(s)-X_{t,Y}(s)\big|^2\Big]ds\\
\leq&\;\E\big[|X-Y|^2\big]+\big(4L^2+4L+8K_1\big)\int_t^{t+\theta}\E\Big[\big|X_{t,X}(s)-X_{t,Y}(s)\big|^2\Big]ds,
\end{align*}
where $K_1=\int_{\msE}\rho^2(s)\lambda(de)\vee\lambda(\msE)$.
Then by the Gronwall lemma \cite{kp}, we obtain
\begin{equation}\label{e5}
\E\big[|X_{t,X}(t+\theta)-X_{t,Y}(t+\theta)|^2\big]\leq  e^{4(L^2+L+2K_1)h} \E\big[|X-Y|^2\big],
\end{equation}
which leads to the inequality \eqref{e3}.
	
By the definition of $Z$, we have
\begin{equation}\begin{aligned}\label{zz}
Z=&\;\int_t^{t+h}\!\Big(b^{X_{t,X}}\big(s,X_{t,X}(s)\big)-b^{X_{t,Y}}\big(s,X_{t,Y}(s)\big)\Big)ds\\
&\;+\int_t^{t+h}\!\Big(\sigma^{X_{t,X}}\big(s,X_{t,X}(s)\big)-\sigma^{X_{t,Y}}\big(s,X_{t,Y}(s)\big)\Big)dW_s\\
&\;+\int_{t}^{t+h}\!\!\!\int_{\msE}\Big(c^{X_{t,X}}\big(s,X_{t,X}(s-),e\big)-c^{X_{t,Y}}\big(s,X_{t,Y}(s-),e\big)\Big)\mu(de,ds)\\
=&\;\int_t^{t+h}\!\Big(\tilde{b}^{X_{t,X}}\big(s,X_{t,X}(s)\big)-\tilde{b}^{X_{t,Y}}\big(s,X_{t,Y}(s)\big)\Big)ds\\
&\;+\int_t^{t+h}\!\Big(\sigma^{X_{t,X}}\big(s,X_{t,X}(s)\big)-\sigma^{X_{t,Y}}\big(s,X_{t,Y}(s)\big)\Big)dW_s\\
&\;+\int_{t}^{t+h}\!\!\!\int_{\msE}\Big(c^{X_{t,X}}\big(s,X_{t,X}(s-),e\big)-c^{X_{t,Y}}\big(s,X_{t,Y}(s-),e\big)\Big)\tilde{\mu}(de,ds).
\end{aligned}\end{equation}
Taking square on both sides of \eqref{zz}
and taking $\mathbb E[\cdot]$ on the derived equation, we get
\begin{align}\label{zzz}
\E[|Z|^2]\leq&\;3\E\Bigg[\Big|\int_t^{t+h}\Big(\tilde{b}^{X_{t,X}}\big(s,X_{t,X}(s)\big)-\tilde{b}^{X_{t,Y}}\big(s,X_{t,Y}(s)\big)\Big)ds\Big|^2\Bigg]\nonumber\\
&\;+3\E\Bigg[\int_t^{t+h}\big|\sigma^{X_{t,X}}\big(s,X_{t,X}(s)\big)-\sigma^{X_{t,Y}}\big(s,X_{t,Y}(s)\big)\big|^2ds\Bigg]\nonumber\\
&\;+3\E\Bigg[\int_t^{t+h}\int_{\msE}\big|c^{X_{t,X}}\left(s,X_{t,X}(s),e\right)-c^{X_{t,Y}}\left(s,X_{t,Y}(s),e\right)\big|^2\lambda(de)ds\Bigg].
\end{align}
Then by using \eqref{e5} and \eqref{zzz}, we deduce
\begin{equation*}\label{zzzz}
\begin{aligned}
\E[|Z|^2]\leq& \;12(L^2\!+\!K_1^2)h\E\Bigg[\int_t^{t+h}\!\!\!\Big(\E\big[|X_{t,X}(s)-X_{t,Y}(s)|^2\big]+|X_{t,X}(s)-X_{t,Y}(s)|^2\Big)ds\Bigg]\\
&\;+6L^2\E\Bigg[\int_t^{t+h}\Big(\E\big[|X_{t,X}(s)-X_{t,Y}(s)|^2\big]+|X_{t,X}(s)-X_{t,Y}(s)|^2\Big)ds\Bigg]\\
&\;+6K_1\E\Bigg[\int_t^{t+h}\Big(\E\big[|X_{t,X}(s)-X_{t,Y}(s)|^2\big]+|X_{t,X}(s)-X_{t,Y}(s)|^2\Big)ds\Bigg]\\
\leq&\; 12(2L^2+2K_1^2+K_1)(1+h)\Bigg[\int_t^{t+h}\E\big[|X_{t,X}(s)-X_{t,Y}(s)|^2\big]ds\Bigg]\\
\leq&\; 12(2L^2+2K_1^2+K_1)(1+h)\Bigg[\int_t^{t+h}\E\big[|X-Y|^2\big]\times e^{4(L^2+L+2K_1)h}ds\Bigg]\\
\leq&\; 12(2L^2+2K_1^2+K_1)(1+h)e^{4(L^2+L+2K_1)h}\E\big[|X-Y|^2\big]h,
\end{aligned}
\end{equation*}
which proves \eqref{e4}. The proof ends.
\end{proof}

\begin{lem}\label{l2}  Under the assumptions ${\rm (A1)} - {\rm (A4)}$, for $k=0,\cdots,N-1$, we have
\begin{equation}\label{lem}
\E\Big[\big|\E\big[X_{t_k,\bar{X}_k}(t_{k+1})-\bar{X}_k\big|\mathcal{F}_{t_k}\big]\big|^2\Big]
\leq C\big(1+\E\big[|\bar{X}_k|^2\big]\big)h^2,
\end{equation}
where $C$ is a positive constant depending on $\lambda(\msE)$,
the function $\rho(e)$  in ${\rm (A3)}$
and the linear growth constant $K$ in ${\rm (A4)}$.
\end{lem}
\begin{proof} By \eqref{ree}, we get
\begin{align*}
\E\big[X_{t_k,\bar{X}_k}(t_{k+1})-\bar{X}_k|\mathcal{F}_{t_k}\big]
=&\;\E\Bigg[\int_{t_k}^{t_{k+1}}b^{X_{t_k,\bar{X}_k}}\big(s,X_{t_k,\bar{X}_k}(s)\big)ds\Big|\mathcal{F}_{t_k}\Bigg]\\
&\;+\E\Bigg[\int_{t_k}^{t_{k+1}}\!\!\!\int_{\msE}c^{X_{t_k,\bar{X}_k}}\big(s,X_{t_k,\bar{X}_k}(s),e\big)\lambda(de)ds\Big|\mathcal{F}_{t_k}\Bigg]\\
=&\;\E\Bigg[\int_{t_k}^{t_{k+1}}\tilde{b}^{X_{t_k,\bar{X}_k}}\big(s,X_{t_k,\bar{X}_k}(s)\big)ds\Big|\mathcal{F}_{t_k}\Bigg].
\end{align*}
Then
\begin{equation}\label{eq:con}
\begin{aligned}
&\;\E\Big[\big|\E\big[X_{t_k,\bar{X}_k}(t_{k+1})-\bar{X}_k|\mathcal{F}_{t_k}\big]\big|^2\Big]\\
\leq &\;\E\bigg[\E\bigg[\Big|\int_{t_k}^{t_{k+1}}\tilde{b}^{X_{t_k,\bar{X}_k}}\big(s,X_{t_k,\bar{X}_k}(s)\big)ds\Big|^2\Big|\mathcal{F}_{t_k}\bigg]\bigg]\\
\leq &\;h\int_{t_k}^{t_{k+1}}\E\bigg[\E\Big[\big|\tilde{b}\big(s,X_{t_k,\bar{X}_k}(s),x\big)\big|^2\Big]\Big|_{x=X_{t_k,\bar{X}_k}(s)}\bigg]ds.
\end{aligned}
\end{equation}
Using the assumption ${\rm (A4)}$, it is easy to have
\begin{align}\label{eq:tildeb}
\big|\tilde{b}(s,X_{t_k,\bar{X}_k}(s),x)\big|^2\leq C\big(1+|X_{t_k,\bar{X}_k}(s)|^2+|x|^2\big),
\end{align}
where $C$ depends on $\lambda(\msE)$, $\rho(e)$ and  $K$.
Then by Theorem \ref{th0} and \eqref{eq:tildeb}, we deduce
\begin{equation}\begin{aligned}\label{eq:estb}
&\;\E\bigg[\E\Big[\big|\tilde{b}\big(s,X_{t_k,\bar{X}_k}(s),x\big)\big|^2\Big]\Big|_{x=X_{t_k,\bar{X}_k}(s)}\bigg]\\
\leq &\;C\Big(1+2\E\big[|X_{t_k,\bar{X}_k}(s)|^2\big]\Big)
\leq C\Big(1+\E\big[|\bar X_k|^2\big]\Big).
\end{aligned}\end{equation}
Combining with \eqref{eq:con} and \eqref{eq:estb}, we obtain
\begin{align*}
&\;\E\Big[\big|\E\big[X_{t_k,\bar{X}_k}(t_{k+1})-\bar{X}_k|\mathcal{F}_{t_k}\big]\big|^2\Big]\\
\leq&\; h\int_{t_k}^{t_{k+1}}C\big(1+\E\big[|\bar{X}_k|^2\big]\big)ds
\leq C\big(1+\E\big[|\bar{X}_k|^2\big]\big)h^2,
\end{align*}
which completes the proof.
\end{proof}




Now, based on the Lemmas \ref{l1} and \ref{l2},
we give the general error estimate theorem for the one-step scheme \eqref{onestep}.
\begin{thm}\label{th51}
	Let $X_{t,X}(t+h)$ be defined as \eqref{ree}.
	If $\bar{X}_{t,X}(t+h)$ satisfies
\begin{subequations}
	\begin{align}
	&\big|\E\big[X_{t,X}(t+h)-\bar{X}_{t,X}(t+h)\big|\mathcal{F}_t\big]\big|
	\leq C^{*} \big(1+\E\big[|X|^2\big]+|X|^2\big)^{\frac{1}{2}}h^{p_1},\label{e1}\\
	&\big(\E\big[|X_{t,X}(t+h)-\bar{X}_{t,X}(t+h)|^2\big|\mathcal{F}_t\big]\big)^\frac{1}{2}
	\leq C^{*} \big(1+\E\big[|X|^2\big]+|X|^2\big)^{\frac{1}{2}}h^{p_2},\label{e2}
	\end{align}
\end{subequations}
where $t\in [0,T-h]$, $p_1$ and $p_2$ are parameters satisfying
$p_2\ge\frac{1}{2}$ and $p_1\ge p_2+\frac{1}{2}$,
and $C^{*}>0$ is a constant independent of
$h$,  $X_{t,X}(t+h)$ and $\bar{X}_{t,X}(t+h)$.
Then for $k=1,\cdots,N$, it holds that
\begin{equation}\label{gerror}
\big(\E\big[|X(t_k)-\bar{X}_{k}|^2\big]\big)^\frac{1}{2}
\leq C\big(1+\E\big[|X_0|^2\big]\big)^{\frac{1}{2}}h^{p_2-\frac{1}{2}},
\end{equation}
where
$C$ is a constant independent of  $h$,  $X_{t,X}(t+h)$ and $\bar{X}_{t,X}(t+h)$.
\end{thm}

\begin{proof}
By Theorem \ref{th0}, Lemma \ref{l2} and the discrete Gronwall lemma \cite{mt},
it is easy to prove that for all $k=0,\cdots,N$,
\begin{equation}\label{momentest}
\E\big[|\bar{X}_k|^2\big]\leq C\big(1+\E\big[|X_0|^2\big]\big).
\end{equation}
%

Then based on Lemmas \ref{l1} and \ref{l2} and the inequality \eqref{momentest},
the proof of Theorem \ref{th51} is similar to that of Theorem 4.1 in \cite{SYZ2017}.
So we omit it here.
\end{proof}

\begin{rem}
Theorem \ref{th51} implies that when the weak local error estimate
of the one-step scheme \eqref{onestep} is of order $p_1$ and	
its strong local error estimate
is of order $p_2$, then
the global strong order of the scheme \eqref{onestep} is $p_2-\frac{1}{2}$.
\end{rem}

\subsection{The error estimates for strong Taylor schemes}\label{c6s44}

In this subsection, utilizing  Theorem \ref{th51}, we prove the
error estimates of Schemes \ref{SSG} and  \ref{SSGcom}
to reveal the orders of strong convergence of strong Taylor schemes.

Let $W_t^0=t$ and
%
$\alpha=(i_1,i_2,\cdots,i_k)\in\mathcal{M}$ be a given multi-index.
Then we have the following two lemmas.
\begin{lem}\label{lem11}
Let  $f^\beta$ and $\beta$ be defined by  \eqref{itf} and \eqref{itpppj}, respectively.
Assume that $\tilde{f}_\alpha^\beta$ and $\tilde{I}_\alpha\big[\tilde{f}^\beta(\cdot)\big]_{t,t+h}$ exist and
$
\big|\tilde{f}_\alpha^\beta(t,x)\big|\leq C\big(1+\E\big[|\beta_t|^2\big]+|x|^2\big)^{1/2}.
$
Then
\begin{equation}\label{ff1}
\E\Big[\Big(\tilde{I}_{\alpha}\big[\tilde{f}_\alpha^{X_{t,X}}\big(\cdot,X_{t,X}(\cdot)\big)\big]_{t,t+h}\Big)^{2}\Big|\F_t\Big]\leq CMh^{l(\alpha)+n(\alpha)},\qquad\qquad\quad\;\;
\end{equation}
and
\begin{equation}\label{ff1w}\begin{aligned}
\Big|\E \Big[\tilde{I}_{\alpha}\big[\tilde{f}_\alpha^{X_{t,X}}\big(\cdot,X_{t,X}(\cdot)\big)\big]_{t,t+h}\big|\mtF_t\Big]\Big|
\left\{\begin{aligned}=&\; 0,& \text{if}~~l(\alpha)\neq n(\alpha),\\
\leq &\; CMh^{l(\alpha)},  &\text{if}~~ l(\alpha)= n(\alpha),
\end{aligned}\right.
\end{aligned}\end{equation}
where $M=\left(1+\E\big[|X|^2\big]+|X|^2\right)^{1/2}$.
\end{lem}
\begin{proof}
If $\alpha=v$, i.e.,  $l(\alpha)+n(\alpha)=0$, we get
\begin{equation*}\begin{aligned}\label{ff}
&\;\E\Big[\Big(\tilde{I}_{\alpha}\big[\tilde{f}_\alpha^{X_{t,X}}\big(\cdot,X_{t,X}(\cdot)\big)\big]_{t,t+h}\Big)^{2}\Big|\F_t\Big]\\
=&\;\E\Big[\big|\tilde{f}^{X_{t,X}}\big(t+h,X_{t,X}(t+h)\big)\big|^2\Big|\F_t\Big]\\
\leq &\;
C\E\Big[1+\E\big[|X_{t,X}(t+h)|^2\big]+|X_{t,X}(t+h)|^2\Big|\F_t\Big]\\
\leq &\; C\big(1+\E[|X|^2]+|X|^2\big),
\end{aligned}\end{equation*}
which leads to \eqref{ff1} with $p(\alpha)=l(\alpha)+n(\alpha)=0$.

Now we consider $l(\alpha)>0$. If $i_k\neq0$, by the It\^o's isometry formula, we have
\begin{equation}\begin{aligned}\label{ref1}
&\;\E\Big[\Big(\tilde{I}_{\alpha}\big[\tilde{f}_\alpha^{X_{t,X}}\big(\cdot,X_{t,X}(\cdot)\big)\big]_{t,t+h}\Big)^{2}\Big|\F_t\Big]\\
=&\;\left\{\begin{aligned}&\int_{t}^{t+h}\E\Big[\Big(\tilde{I}_{\alpha-}
\big[\tilde{f}_\alpha^{X_{t,X}}\big(\cdot,X_{t,X}(\cdot)\big)\big]_{t,s}\Big)^{2}\Big|\F_t\Big]ds,
&\text{if}~~ i_k=1,\\
&\int_{t}^{t+h}\!\!\!\int_{\msE}\E\Big[\Big(\tilde{I}_{\alpha-}\big[\tilde{f}_\alpha^{X_{t,X}}
\big(\cdot,X_{t,X}(\cdot)\big)\big]_{t,s}\Big)^{2}\Big|\F_t\Big]\lambda(de)ds,
&\text{if}~~ i_k=-1.
\end{aligned}\right.
\end{aligned}\end{equation}
If $i_k=0$, by the Holder's inequality, we obtain
\begin{equation}\begin{aligned}\label{ref2}
& \;\E\Big[\Big(\tilde{I}_{\alpha}\big[\tilde{f}_\alpha^{X_{t,X}}\big(\cdot,X_{t,X}(\cdot)\big)\big]_{t,t+h}\Big)^{2}\Big|\F_t\Big]\\
=& \; \E\bigg[\Big|\int_{t}^{t+h}\tilde{I}_{\alpha-}\big[\tilde{f}_\alpha^{X_{t,X}}\big(\cdot,X_{t,X}(\cdot)\big)\big]_{t,s}ds\Big|^2\Big|\F_t\bigg]\\
\leq &\; h\int_{t}^{t+h}\E\Big[\Big(\tilde{I}_{\alpha-}\big[\tilde{f}_\alpha^{X_{t,X}}\big(\cdot,X_{t,X}(\cdot)\big)\big]_{t,s}\Big)^{2}\Big|\F_t\Big]ds.
\end{aligned}\end{equation}
Then combining with \eqref{ref1} and \eqref{ref2}, we deduce the recurrence relation
\begin{equation*}\begin{aligned}
p(\alpha)= p(\alpha-)+\big(1+\mathbb{I}_{\{i_{k}=0\}}\big)
=\sum\limits_{j=1}^{k}\big(1+\mathbb{I}_{\{i_j=0\}}\big)=l(\alpha)+n(\alpha),
\end{aligned}\end{equation*}
which implies that \eqref{ff1} holds ture.
%

Similarly, we can prove \eqref{ff1w}.
The proof ends.
\end{proof}

\begin{lem}\label{lem33}
Let $f^\beta$ and $\beta$ be defined by \eqref{itf} and \eqref{itpppj}, respectively.
Assume that $f_\alpha^\beta$ and $I_\alpha\big[f^\beta(\cdot)\big]_{t,t+h}$ exist and
$
\big|f_\alpha^\beta(t,x)\big|\leq C\big(1+\E\big[|\beta_t|^2\big]+|x|^2\big)^{1/2}.
$
Then
\begin{equation}\label{ff3}
\begin{aligned}
\E\Big[\Big(I_{\alpha}\big[f_\alpha^{X_{t,X}}\big(\cdot,X_{t,X}(\cdot)\big)\big]_{t,t+h}\Big)^{2}\Big|\F_t\Big]\leq CMh^{l(\alpha)+n(\alpha)},
\qquad\qquad\qquad\;\;\;\;\;\;\;\;
\end{aligned}\end{equation}
and
\begin{equation}\label{ff4}\begin{array}{lll}
 \Big|\E \Big[I_{\alpha}\big[f_\alpha^{X_{t,X}}\big(\cdot,X_{t,X}(\cdot)\big)\big]_{t,t+h}\big|\mtF_t\Big]\Big|
\left\{\begin{aligned}
=&\; 0,& \text{if}~l(\alpha)\neq n(\alpha)+s(\alpha),\\
\leq &\; CMh^{l(\alpha)},  &\text{if}~ l(\alpha)= n(\alpha)+s(\alpha),
\end{aligned}\right.\end{array}
\end{equation}
where $M=\big(1+\E\big[|X|^2\big]+|X|^2\big)^{1/2}$.
\end{lem}
\begin{proof}
By Lemma \ref{lem11} and the relationship $\tilde{\mu}(de,dt)=\mu(de,dt)-\lambda(de)dt$,
it is easy to prove \eqref{ff3} and \eqref{ff4}. The proof ends.
\end{proof}

Based on Theorem \ref{th51} and Lemma \ref{lem33},
we prove the error estimate of the strong order $\gamma$ Taylor scheme
 in the following theorem.

\begin{thm}\label{strong:thm}
Let $X(t)$ and $\bar X_{k}$ be the solutions of the MSDEJ \eqref{MSDEJeq}
and the strong order $\gamma$ Taylor scheme \ref{SSG}, respectively.
Let $f(t,x',x)=x$ and
assume that $f^{X_{t,X}}\big(s,X_{t,X}(s)\big)$
has the It\^o-Taylor expansion \eqref{ite} with $\mathcal{A}=\mathcal{A}_\gamma$ and
$$\big|f_\alpha^\beta(t,x)\big|\leq C\big(1+\E\big[|\beta_t|^2\big]+|x|^2\big)^{1/2},\;\;\;(t,x) \in [0,T]\times\R^{d}$$
for all $\alpha\in \mathcal{A}_\gamma\cup\mathcal B(\mathcal A_\gamma)$
with $\beta$ defined by \eqref{itpppj}. Then it holds that
$$\max_{k\in \{1,2,\cdots,N\}}\mathbb E\big[|X_{t_{k}}-\bar X_{k}|^{2}\big]
\le C\big(1+\mathbb E\big[|X_{0}|^{2}\big]\big) \big(\Delta t\big)^{2\gamma},$$
\end{thm}
where $\Delta t=\Delta t_k=T/N$ for $k=0,1,\dots,N-1$.
\begin{proof}
By the It\^o-Taylor expansion \eqref{ite},  we have
\begin{equation}\label{e7000}
X_{t,X}(t+h)=X+\sum_{\alpha\in\mathcal{A_\gamma}\backslash v}I_\alpha\left[f_\alpha^{X_{t,X}}(t,X)\right]_{t,t+h}
+R^{\gamma},
\end{equation}
where \[R^{\gamma}=\sum\limits_{\alpha\in\mathcal{B(A_\gamma)}}I_\alpha\left[f_\alpha^{X_{t,X}}\big(\cdot,X_{t,X}(\cdot)\big)\right]_{t,t+h}.\]
Moreover, the strong order $\gamma$ Taylor scheme \ref{SSG} can be written as
\begin{equation}\label{e700}
\bar{X}_{t,X}(t+h)=X+\sum_{\alpha\in\mathcal{A_\gamma}\backslash v}I_\alpha\left[f_\alpha^{X_{t,X}}(t,X)\right]_{t,t+h}.
\end{equation}
Then we subtract  \eqref{e700} from \eqref{e7000} and obtain
\begin{equation}\label{e888}
\begin{aligned}
R^\gamma = X_{t,X}(t+h)-\bar{X}_{t,X}(t+h).
\end{aligned}
\end{equation}
According to Lemma \ref{lem33}, we deduce
\begin{align*}
\E\Big[|R^\gamma|^2\big|\F_t\Big]&=\E\bigg[\Big|\sum\limits_{\alpha\in\mathcal{B(A_\gamma)}}I_\alpha
\big[f_\alpha^{X_{t,X}}\big(\cdot,X_{t,X}(\cdot)\big)\big]_{t,t+h}\Big|^2\Big|\mtF_t\bigg]\\
&\leq C\sum\limits_{\alpha\in\mathcal{B(A_\gamma)}}\E\Big[\big|I_\alpha\big[f_\alpha^{X_{t,X}}\big(\cdot,X_{t,X}(\cdot)\big)\big]_{t,t+h}\big|^2\Big|\F_t\Big]\\
&\leq C\sum\limits_{\alpha\in\mathcal{B(A_\gamma)}} \big(1+\E\big[|X|^2\big]+|X|^2\big)h^{l(\alpha)+n(\alpha)}\\
& \le C\big(1+\E\big[|X|^2\big]+|X|^2\big)h^{2p_2},
\end{align*}
where $2p_2=\min\limits_{\alpha\in\mathcal{B(A_\gamma)}}\{l(\alpha)+n(\alpha)\}$.
Since
$\mathcal{A}_\gamma = \big\{\alpha \in\mathcal{M}: l(\alpha) + n(\alpha)
\leq 2\gamma ~\text{or}~  l(\alpha) = n(\alpha) = \gamma + \frac{1}{2}\big\}$,
then we get
$$p_2=\frac12 \min\limits_{\alpha\in\mathcal{B(A_\gamma)}}\{l(\alpha)+n(\alpha)\}=\gamma+\frac{1}{2}.$$
	
Now we prove $p_1\geq p_2+\frac{1}{2}$.
By Lemma \ref{lem33}, we can deduce
\begin{align*}
\Big|\E\big[R^\gamma\big|\F_t\big]\Big|&=\Big|\E\Big[\sum\limits_{\alpha\in\mathcal{B(A_\gamma)}}
I_\alpha\big[f_\alpha^{X_{t,X}}\big(\cdot,X_{t,X}(\cdot)\big)\big]_{t,t+h}\big|\F_t\Big]\Big|\\
&\leq\sum\limits_{\alpha\in\mathcal{B(A_\gamma)}}\Big|\E\Big[I_\alpha\big[f_\alpha^{X_{t,X}}\big(\cdot,X_{t,X}(\cdot)\big)\big]_{t,t+h}\big|\F_t\Big]\Big|\\
&=\sum\limits_{\tiny{\begin{array}{c}
\alpha\in\mathcal{B(A_\gamma)}\\l(\alpha)=n(\alpha)+s(\alpha)
\end{array}}}\Big|\E\Big[I_\alpha\big[f_\alpha^{X_{t,X}}\big(\cdot,X_{t,X}(\cdot)\big)\big]_{t,t+h}\big|\F_t\Big]\Big|\\
&\leq C\sum\limits_{\tiny{\begin{array}{c}
\alpha\in\mathcal{B(A_\gamma)}\\l(\alpha)=n(\alpha)+s(\alpha)
\end{array}}}\big(1+\E\big[|X|^2\big]+|X|^2\big)^{1/2}h^{l(\alpha)}\\
& \le  C\big(1+\E\big[|X|^2\big]+|X|^2\big)^{1/2}
h^{p_1},
\end{align*}
where $p_1=\min\limits_{\alpha\in\mathcal{B(A_\gamma)}}\{l(\alpha):l(\alpha)=n(\alpha)+s(\alpha)\}$.
Simple calculation yields
\begin{equation*}
p_1=\min\limits_{\alpha\in\mathcal{B(A_\gamma)}}\{l(\alpha):l(\alpha)=n(\alpha)+s(\alpha)\}=	
\begin{cases}
\gamma+1,\; \gamma=1,2,\cdots\\
\gamma+\frac{3}{2},\; \gamma=0.5,1.5,\cdots,
\end{cases}
\end{equation*}
which implies that $p_1\geq p_2+\frac{1}{2}$.
Then by Theorem \ref{th51}, we complete the proof.
\end{proof}

From Theorem \ref{th51}, we come to the conclusion that the order of strong convergence of
the strong order $\gamma$ Taylor scheme \ref{SSG} is $\gamma$.
Moreover, by using Theorem \ref{th51} and Lemma \ref{lem11},
we can prove that
the order of strong convergence of the
compensated strong order $\gamma$ Taylor scheme \ref{SSGcom} is also $\gamma$.

\begin{coro}\label{coro}
Let $X(t)$ and $\bar X_{k}$ be the solutions of the MSDEJ \eqref{MSDEJeq}
and the compensated strong order $\gamma$  Taylor scheme \ref{SSGcom}, respectively. 
Let $f(t,x',x)=x$ and assume that $f^{X_{t,X}}\big(s,X_{t,X}(s)\big)$
has the It\^o-Taylor expansion \eqref{itecompensated} with $\mathcal{A}=\mathcal{A}_\gamma$ and
$$\big|\tilde{f}_\alpha^\beta(t,x)\big|\leq C\big(1+\E\big[|\beta_t|^2\big]+|x|^2\big)^{1/2},\;\;\;(t,x) \in [0,T]\times\R^{d}$$
for all $\alpha\in \mathcal{A}_\gamma\cup\mathcal B(\mathcal A_\gamma)$
with $\beta$ defined by \eqref{itpppj}. Then it holds that
$$\max_{k\in \{1,2,\cdots,N\}}\mathbb E\big[|X_{t_{k}}-\bar X_{k}|^{2}\big]
\le C\big(1+\mathbb E\big[|X_{0}|^{2}\big]\big) \big(\Delta t\big)^{2\gamma}.$$
\end{coro}

\section{Error estimates for weak Taylor schemes}\label{sec:weakerr}
In this section, we focus on the
error estimates of the weak order $\eta$ It\^o-Taylor Scheme \ref{SWG} and
the compensated weak order $\eta$ It\^o-Taylor Scheme \ref{SWGcom}.
For this purpose, we first present
some useful lemmas as below.


Let $C_P^{k,2k,2k}([0,T]\times\R^d\times\R^d;\R)$ be
the set of functions $\varphi(t,x',x):[0,T]\times\R^d\times\R^d\rightarrow\R$ such that all their
derivatives with respect to $t$, $x'$ and $x$ up to $k$, $2k$ and $2k$, respectively,
are continuous and of polynomial growth.

For a given $\eta\in\{1,2,\dots\}$ and function $g\in C_P^{2\eta}(\R^d;\R)$,
we define
\begin{equation}\label{def:komgo}
u(s,y)=\E\big[g\big(X_T^{s,y}\big)\big],
\end{equation}
where $(s,y)\in[0,T]\times\R^d$ and $X_t^{s,y},\;s\le t\le T$, is the solution of the MSDEJ \eqref{MSDEJeq}
starting from $(s,y)$.
Then we get
\begin{equation}\label{komogoeq}
u(0,X_0)=\E\big[g\big(X_T^{0,X_0}\big)\big]=\E\big[g(X_T)\big].
\end{equation}

Then we have the following Kolmogorov backward  equation for MSDEJs.

\begin{lem}[Kolmogorov backward  equation]\label{weak:lem1}
Assume that the
coefficients of
the MSDEJ
\eqref{MSDEJeq} have the components
$b^k,\sigma^{k,j},c^{k}\in C_P^{\eta,2\eta,2\eta}\big([0,T]\times\R^d\times\R^d;\R\big)$
for $1\le k\leq d$ and $1\le j\le m$
with uniformly bounded derivatives.
Then the functional $u$
define in \eqref{def:komgo}
is the unique solution
of the nonlocal Kolmogorov backward partial integral differential equations
\begin{equation}\label{komogorov}
\left\{\begin{aligned}
&\tilde{L}^0u(s,y)=0, &(s,y)\in[0,T)\times\R^d,\\
&u(T,y)=g(y), &y\in\R^d,
\end{aligned}
\right.
\end{equation}
where $\tilde{L}^0$ is defined by \eqref{tildel0}.
Moreover, we have
\begin{equation}\label{komogorov:smooth}
u(s,\cdot)\in C_P^{2\eta}(\R^d;\R),\;\;\;0\le s\le T.
\end{equation}
\end{lem}

\begin{proof}
The proof of Lemma \ref{weak:lem1} is similar to that of
Lemma 12.3.1 in \cite{PB2010}.
So we omit it here.
The readers are referred to \cite{PB2010} for more details.
\end{proof}
Based on the Kolmogorov backward  equation,
we shall prove the error estimates of Schemes \ref{SWG}
and \ref{SWGcom} to reveal the orders of weak convergence of weak Taylor schemes.
To proceed, we introduce the following two lemmas.

For a given $x\in\R$ and $p\in \mathbb{N}$, we denote by $[x]$ the integer part of $x$
and $\mathcal{A}_p$ the set of multi-indices $\alpha = (j_1,\dots,j_l)$
of length $l \le p$ with
components $j_i\in\{-1,0\}$, for $i\in\{1,\dots,l\}$.

\begin{lem}\label{lemmaw1}
Let $X_t$ be the solution of the MSDEJ \eqref{MSDEJeq}, and
$\rho$ and $\tau$ be two stopping times with $\tau$ being $\F_{\rho}$-measurable
and $0\le\rho\le\tau\le T$ a.s..
Given $\alpha\in\mathcal{M}$, let $p=l(\alpha)-\big[\frac{l(\alpha)+n(\alpha)}{2}\big]$
and
$f(t,x)\in C^{p,2p}([\rho,\tau]\times\R^d;\R)$ be a $\F_t$-adapted process
such that for any $\alpha\in \mathcal{A}_p$
$$\E\big[˜\big(f_{\alpha}(t,X_t)\big)^2\big|\F_{\rho}\big]\le K,\;\;\;a.s.,\;\;\;t\in[\rho,\tau],$$
for some constant $K$.
Moreover, for an adapted process $g(\cdot)=g(\cdot,e)$
with $e\in \msE^{s(\alpha)}$,
if $\E\big[g(t,e)^2\big|\F_{\rho}\big]< +\infty$ a.s. for
$t\in[\rho,\tau]$, then
\begin{equation}\label{ineq1}
\Big|\E\Big[f(\tau,X_{\tau})I_{\alpha}[g(\cdot)]_{\rho,\tau}\big|\F_{\rho}\Big]\Big|
\le C_1(\tau-\rho)^{l(\alpha)},
\end{equation}
and
\begin{equation}\label{ineq111}
\Big|\E\Big[f(\tau,X_{\tau})\tilde{I}_{\alpha}[g(\cdot)]_{\rho,\tau}\big|\F_{\rho}\Big]\Big|
\le C_2(\tau-\rho)^{l(\alpha)},
\end{equation}
where the positive constants $C_1$ and $C_2$ do not depend on $(\tau-\rho)$.
\end{lem}
%
%
\begin{proof}
It is obvious that \eqref{ineq1} holds for $|\alpha|=0$. Suppose that
\eqref{ineq1} holds for all $|\alpha|\le l$.
Now we consider $\alpha=(j_1,\dots,j_{l+1})$ with $j_{l+1}=-1$ and obtain
\begin{equation}\label{gl}
I_{\alpha}[g(\cdot)]_{\rho,\tau}=\int_{\rho}^{\tau}\!\!\int_{\msE}I_{\alpha-}[g(\cdot)]_{\rho,s-}\mu(de,ds).
\end{equation}
Then by the It\^o's formula \eqref{mitj}, we get
\begin{align}\label{ineq2}
f(\tau,X_{\tau})&I_{\alpha}[g(\cdot)]_{\rho,\tau}=
\int_{\rho}^{\tau}L^0f(s,X_s)I_{\alpha}[g(\cdot)]_{\rho,s}ds
+\int_{\rho}^{\tau}\overrightarrow{L}^1f(s,X_s)I_{\alpha}[g(\cdot)]_{\rho,s}dW_s\nonumber\\
&\qquad\qquad\quad\;+\int_{\rho}^{\tau}\!\!\int_{\msE}\Big(f(s,X_{s})I_{\alpha}[g(\cdot)]_{\rho,s}
-f(s,X_{s-})I_{\alpha}[g(\cdot)]_{\rho,s-}\Big)\mu(de,ds).
\end{align}
When $s$ is a jump time, we have
\begin{equation*}\label{fact1}\begin{aligned}
f(s,X_{s})=&\;L_e^{-1}f(s,X_{s-})+f(s,X_{s-}),\\
I_{\alpha}[g(\cdot)]_{\rho,s}=&\;I_{\alpha}[g(\cdot)]_{\rho,s-}+I_{\alpha-}[g(\cdot)]_{\rho,s-}.
\end{aligned}\end{equation*}
By inserting the above equations into \eqref{ineq2}, we deduce
\begin{align}\label{ineq3}
f(\tau,X_{\tau})&I_{\alpha}[g(\cdot)]_{\rho,\tau}=
\int_{\rho}^{\tau}L^0f(s,X_s)I_{\alpha}[g(\cdot)]_{\rho,s}ds
+\int_{\rho}^{\tau}\overrightarrow{L}^1f(s,X_s)I_{\alpha}[g(\cdot)]_{\rho,s}dW_s\nonumber\\
&\qquad+\int_{\rho}^{\tau}\!\!\int_{\msE}\Big(f(s,X_{s})I_{\alpha-}[g(\cdot)]_{\rho,s-}
+L_{e}^{-1}f(s,X_{s-})I_{\alpha}[g(\cdot)]_{\rho,s-}\Big)\mu(de,ds),
\end{align}
which leads to
\begin{align}\label{ineq4}
\Big|\E\Big[f(\tau,X_{\tau})I_{\alpha}[g(\cdot)]_{\rho,\tau}\big|\F_{\rho}\Big]\Big|
=&\;\Big|\int_{\rho}^{\tau}\E\Big[L^0f(s,X_s)I_{\alpha}[g(\cdot)]_{\rho,s}\big|\F_{\rho}\Big]ds\nonumber\\
&\;+\int_{\rho}^{\tau}\!\!\int_{\msE}\E\Big[f(s,X_{s})I_{\alpha-}[g(\cdot)]_{\rho,s}\big|\F_{\rho}\Big]\lambda(de)ds\nonumber\\
&\;+\int_{\rho}^{\tau}\!\!\int_{\msE}\E\Big[L_{e}^{-1}f(s,X_{s})I_{\alpha}[g(\cdot)]_{\rho,s}\big|\F_{\rho}\Big]\lambda(de)ds\Big|\nonumber\\
\le&\;C(\tau-\rho)^{l+1}
+\int_{\rho}^{\tau}\Big|\E\Big[L^0f(s,X_s)I_{\alpha}[g(\cdot)]_{\rho,s}\big|\F_{\rho}\Big]\Big|ds\nonumber\\
&\;+\int_{\rho}^{\tau}\!\!\int_{\msE}\Big|\E\Big[L_{e}^{-1}f(s,X_{s})I_{\alpha}[g(\cdot)]_{\rho,s}\big|\F_{\rho}\Big]\Big|\lambda(de)ds.
\end{align}
Moreover, by the conditions of Lemma \ref{lemmaw1}, we have
\begin{align*}
&L^0f(t,x)\in C^{p-1,2(p-1)},\qquad\;\; L_e^{-1}f(t,x)\in C^{p,2p},\\
&\E\big[˜\big(L^0f_{\alpha}(t,X_t)\big)^2\big|\F_{\rho}\big]\le K,\quad
\E\big[˜\big(L_e^{-1}f_{\alpha}(t,X_t)\big)^2\big|\F_{\rho}\big]\le K,
\end{align*}
and
$$\E\big[\big(I_{(\alpha(1))}[g(\cdot)]_{\rho,t}\big)^2\big|\F_{\rho}\big]< +\infty,\quad \text{for}\quad t\in[\rho,\tau].$$
Then we can repeatedly apply \eqref{ineq4} $p$ times to get
\begin{equation}\label{ineq5}\begin{aligned}
&\;\Big|\E\Big[f(\tau,X_{\tau})I_{\alpha}[g(\cdot)]_{\rho,\tau}\big|\F_{\rho}\Big]\Big|\\
\le&\;\sum_{\beta\in\mathcal{A}_p,l(\beta)=p}\int_{\rho}^{\tau}\cdots\int_{\rho}^{s_2}\int_{\msE}\cdots\int_{\msE}
\Big|\E\Big[f_{\beta}(s_1,X_{s_1})I_{\alpha}[g(\cdot)]_{\rho,s_1}\big|\F_{\rho}\Big]\Big|\\
&\;\;\;\qquad\qquad\qquad\cdot\lambda(de_1)\cdots\lambda(de_{s(\beta)})ds_1\cdots ds_p+C(\tau-\rho)^{l+1}.
\end{aligned}\end{equation}
Using the conditions of Lemma \ref{lemmaw1}, for $s\in[\rho,\tau]$, we deduce
\begin{equation*}
\E\left[\big|I_{\alpha}[g(\cdot)]_{\rho,s}\big|^2\Big|\F_{\rho}\right]
\le C(s-\rho)^{l(\alpha)+n(\alpha)},
\end{equation*}
which implies that
\begin{equation}\label{ineq6}\begin{aligned}
\Big|\E\Big[f_{\beta}(s,X_{s})I_{\alpha}[g(\cdot)]_{\rho,s}\big|\F_{\rho}\Big]\Big|^2
\le&\;\E\left[\big|f_{\beta}(s,X_{s})\big|^2\Big|\F_{\rho}\right]\E\left[\big|I_{\alpha}[g(\cdot)]_{\rho,s}\big|^2\Big|\F_{\rho}\right]\\
\le&\;C(s-\rho)^{l(\alpha)+n(\alpha)}.
\end{aligned}\end{equation}
Then by \eqref{ineq5} and \eqref{ineq6}, we obtain
\begin{equation}\label{ineq7}\begin{aligned}
&\;\Big|\E\Big[f(\tau,X_{\tau})I_{\alpha}[g(\cdot)]_{\rho,\tau}\big|\F_{\rho}\Big]\Big|\\
\le&\;C\int_{\rho}^{\tau}\cdots\int_{\rho}^{s_2}
(s_1-\rho)^{\frac{l(\alpha)+n(\alpha)}{2}}ds_1\cdots ds_p+C(\tau-\rho)^{l+1}\\
\le&\;C(\tau-\rho)^{p+\frac{l(\alpha)+n(\alpha)}{2}}+C(\tau-\rho)^{l+1}.
\end{aligned}\end{equation}
Since $p=l(\alpha)-\big[\frac{l(\alpha)+n(\alpha)}{2}\big]$,
then by \eqref{ineq7}, we have
\begin{equation*}
\Big|\E\Big[f(\tau,X_{\tau})I_{\alpha}[g(\cdot)]_{\rho,\tau}\big|\F_{\rho}\Big]\Big|
\le C(\tau-\rho)^{l+1}.
\end{equation*}

Similarly, we can prove \eqref{ineq1} for $\alpha=(j_1,\dots,j_{l+1})$ with $j_{l+1}=0\;\text{or}\;1$.
Then by using \eqref{ineq1} and the relationship $\tilde{\mu}(de,ds)=\mu(de,ds)-\lambda(de)ds$,
we can get \eqref{ineq111}. The proof ends.
\end{proof}
\begin{lem}\label{lemmaw2}
Let $\rho$ and $\tau$ be two stopping times with $\tau$ being $\F_{\rho}$-measurable
and $0\le\rho\le\tau\le T$, a.s..
Given $\alpha\in\mathcal{M}$
and $\{g(t,e),t\in[\rho,\tau]\}$
with $e\in \msE^{s(\alpha)}$ is an adapted process.
If $g(t,e)$ is $2^{s(\alpha)+3}q$ integrable for a given $q\in \mathbb{N^+}$, then
for any square integrable adapted process $\{h(t),t\in[\rho,\tau]\}$, it holds that
\begin{equation}
\left|\E\left[h(\tau)\big|I_{\alpha}[g(\cdot)]_{\rho,\tau}\big|^{2q}\big|\F_{\rho}\right]\right|
\le C_1(\tau-\rho)^{q\big(l(\alpha)+n(\alpha)-s(\alpha)\big)+s(\alpha)},
\end{equation}
and
\begin{equation}
\left|\E\left[h(\tau)\big|\tilde{I}_{\alpha}[g(\cdot)]_{\rho,\tau}\big|^{2q}\big|\F_{\rho}\right]\right|
\le C_2(\tau-\rho)^{q\big(l(\alpha)+n(\alpha)-s(\alpha)\big)+s(\alpha)},
\end{equation}
where the positive constants $C_1$ and $C_2$ do not depend on $(\tau-\rho)$.
\end{lem}
\begin{proof}
The proof of Lemma \ref{lemmaw2} can be found in Lemma 3.2 in \cite{LL2010}
and Lemma 4.5.5 in \cite{PB2010}.
We omit it here.
\end{proof}

Based on Lemmas \ref{weak:lem1} - \ref{lemmaw2},
we now prove the error estimates of the weak order $\eta$ It\^o-Taylor Scheme
\ref{SWG} in the following theorem.

\begin{thm}\label{weak:thm}
Let $X_t$ and $X_{k}$ be the solutions of the MSDEJ \eqref{MSDEJeq}
and the weak order $\eta$ It\^o-Taylor scheme \ref{SWG}, respectively.
Assume that $\E\big[|X_0|^q\big]<\infty$ for $q\ge 1$ and
$b^{k},\sigma^{k,j},c^k\in C_P^{\eta+1,2(\eta+1),2(\eta+1)}\big([0,T]\times\R^d\times\R^d;\R\big)$
are Lipschitz continuous for $1\le k\le d$ and $1\le j\le m$.
Let the coefficients $f_{\alpha}$ with $f(t,x',x)=x$
satisfy
\begin{equation}\label{cond:linear}
\big|f_{\alpha}^{\beta}(t,x)\big|\le K\big(1+\E\big[|\beta_t|\big]+|x|\big),\qquad (t,x)\in[0,T]\times\R^d
\end{equation}
for all  $\alpha\in\Gamma_{\eta}\bigcup\mathcal{B}(\Gamma_{\eta})$
with $K>0$ being a constant and $\beta_t$ defined by \eqref{itpppj}.
Then for any function $g\in C_P^{2(\eta+1)}(\R^d;\R)$,
it holds that
\begin{equation}\label{weakerror}
\big|\E\big[g(X_T)-g(X_{N})\big]\big|\le C(\Delta t)^\eta,
\end{equation}
where $C$ is a positive constant independent of $\Delta t$.
\end{thm}
\begin{proof}
For simplicity, we consider $d =m= 1$. The proof of the general case is similar.
According to \eqref{komogoeq} and  \eqref{komogorov}, it holds that
\begin{equation}\label{errorh}\begin{aligned}
H=&\;\big|\E\big[g(X_{N})\big]-\E\big[g(X_T)\big]\big|\\
=&\;\big|\E\big[u(T,X_{N})\big]-u(0,X_0)\big|.
\end{aligned}\end{equation}
By the It\^o-Taylor expansion  \eqref{ite} and  \eqref{komogorov},
we deduce
\begin{equation}\label{komocoro}
\E\big[u(t,X_{t}^{s,y})-u(s,y)\big|\F_{s}\big]=0
\end{equation}
for any $0\le s\le t\le T$ and $y\in\R^d$.
Then by \eqref{komogorov:smooth}, \eqref{errorh} and \eqref{komocoro}, we get
\begin{equation}\label{errorh0}\begin{aligned}
H=&\;\Big|\E\Big[\sum_{k=1}^{N}\big(u(t_k,X_{k})-u(t_{k-1},X_{k-1})\big)\Big]\Big|\\
=&\;\Big|\E\Big[\sum_{k=1}^{N}\big(u(t_k,X_{k})-u(t_{k},X_{t_k}^{t_{k-1},X_{k-1}})\big)\Big]\Big|\\
\le&\; H_1+ H_2,
\end{aligned}\end{equation}
where
\begin{equation}\label{errorh_1}
H_1=\Big|\E\Big[\sum_{k=1}^{N}\frac{\partial u}{\partial y}\big(t_k,X_{t_k}^{t_{k-1},X_{k-1}}\big)
\big(X_k-X_{t_k}^{t_{k-1},X_{k-1}}\big)\Big]\Big|,
\end{equation}
\begin{equation}\label{errorh_2}\begin{aligned}
H_2=&\Big|\E\Big[\sum_{k=1}^{N}
\frac12\frac{\partial^2 u}{\partial y^2}\big(t_k,X_{t_k}^{t_{k-1},X_{k-1}}
+\theta_k(X_k-X_{t_k}^{t_{k-1},X_{k-1}})\big)\\
&\qquad\qquad\qquad\times\big(X_k-X_{t_k}^{t_{k-1},X_{k-1}}\big)^2\Big]\Big|
\end{aligned}\end{equation}
with the parameter $\theta_k\in(0,1)$.

By Theorem \ref{th33} and
Scheme \ref{SWG}, we have
\begin{equation}\label{expsch}
X_{t_k}^{t_{k-1},X_{k-1}}-X_k=\sum_{\alpha\in\mathcal{B}(\Gamma_{\eta})}
I_{\alpha}\Big[f_{\alpha}^X\big(\cdot,X_{\cdot}^{t_{k-1},X_{k-1}}\big)\Big]_{t_{k-1},t_k}.
\end{equation}
Utilizing the estimate \eqref{momentest} and the condition \eqref{cond:linear},
one can prove that for every $p\in\{1,2,\dots\}$,
there exist constants $C$ and $r$ such that for every $q\in \{1,\dots,p\}$
\begin{equation}\label{est:sou}
\E\Big[\max_{0\le n\le N}|X_n|^{2q}\Big]\le C\big(1+|X_0|^{2r}\big).
\end{equation}
Since $l(\alpha)=\eta+1$ for $\alpha\in\mathcal{B(\eta)}$,
then we obtain  $2\eta+1\ge 2p$ for $p=l(\alpha)-\big[\frac{l(\alpha)+n(\alpha)}{2}\big]$.
Then by \eqref{komogorov:smooth} and \eqref{est:sou}, for $\alpha\in\mathcal{A}_p$
and $k=1,\dots,N$, we deduce
\begin{align}
\E\Big[˜\big(U_{\alpha}(t_k,X_{t_k})\big)^2\big|\F_{\rho}\Big]<&\;+\infty,\;\;\;a.s.,\label{condition2}\\
\E\Big[˜\big(V(t_k,X_{t_k})\big)^2\big|\F_{\rho}\Big]\;\,<&\;+\infty,\;\;\;a.s.,\label{condition3}
\end{align}
where
\begin{align*}
U(t_k,X_{t_k})=&\;\frac{\partial u}{\partial y}\big(t_k,X_{t_k}^{t_{k-1},X_{k-1}}\big),\\
V(t_k,X_{t_k})=&\;\frac{\partial^2 u}{\partial y^2}\big(t_k,X_{t_k}^{t_{k-1},X_{k-1}}
+\theta_k(X_k-X_{t_k}^{t_{k-1},X_{k-1}})\big).
\end{align*}
Moreover, by Theorem \ref{th0} and the condition \eqref{cond:linear},
 we get
\begin{equation}\label{condition1}
\E\Big[\big|f_{\alpha}^X\big(z,X_{z}^{t_{k-1},X_{k-1}}\big)\big|^2\Big|\F_{t_{k-1}}\Big]
\le C\big(1+\E\big[|X_0|^2\big]+|X_0|^2\big).
\end{equation}
Then based on \eqref{expsch},  \eqref{condition2} and \eqref{condition1},
by Lemma \ref{lemmaw1}, we have
\begin{equation}\label{erresth1}
\begin{aligned}
H_1\le &\;\E\bigg[\sum_{k=1}^{N}\!\sum_{\{\alpha:l(\alpha)=\eta+1\}}\!
\Big|\E\Big[\frac{\partial u}{\partial y}\big(t_k,X_{t_k}^{t_{k-1},X_{k-1}}\big)
I_{\alpha}\Big[f_{\alpha}^X\big(\cdot,X_{\cdot}^{t_{k-1},X_{k-1}}\big)\Big]_{t_{k-1},t_k}\Big|\F_{t_{k-1}}\Big]\Big|\bigg]\\
\le&\;C\E\bigg[\sum_{k=1}^{N}\sum_{\{\alpha:l(\alpha)=\eta+1\}}(t_k-t_{k-1})^{\eta+1}\bigg]
\le C(\Delta t_k)^{\eta}.
\end{aligned}\end{equation}
%
%
%
Similarly, based on \eqref{expsch}, \eqref{condition3} and \eqref{condition1}, by Lemma \ref{lemmaw2} with $q=1$,
we deduce
\begin{equation}\label{erresth2}\begin{aligned}
H_2\le&\;C\E\bigg[\sum_{k=1}^{N}\sum_{\{\alpha:l(\alpha)=\eta+1\}}
\E\Big[\Big|\frac{\partial^2 u}{\partial y^2}\big(t_k,X_{t_k}^{t_{k-1},X_{k-1}}
+\theta_k(X_k-X_{t_k}^{t_{k-1},X_{k-1}})\big)\Big|\\
&\qquad\qquad\qquad\qquad\qquad\times\Big|I_{\alpha}\Big[f_{\alpha}^X\big(\cdot,X_{\cdot}^{t_{k-1},X_{k-1}}\big)\Big]_{t_{k-1},t_k}\Big|^2
\Big|\F_{t_{k-1}}\Big]\bigg]\\
\le&\;C(\Delta t_k)^{\eta}.
\end{aligned}\end{equation}

Then by using  \eqref{errorh}, \eqref{erresth1} and \eqref{erresth2}, we get
 \eqref{weakerror}. The proof ends.
\end{proof}


From Theorem \ref{weak:thm}, we can conclude that the order of weak convergence of
the weak order $\eta$ Taylor scheme \ref{SWG} is $\eta$.
Moreover, based on Lemmas \ref{weak:lem1} - \ref{lemmaw2},
we can prove that the order of weak convergence of the
compensated weak order $\eta$ Taylor scheme \ref{SWGcom} is also $\eta$.

%
\begin{coro}\label{weak:coro}
Let $X_t$ and $X_{k}$ be the solutions of the MSDEJ \eqref{MSDEJeq}
and the compensated weak order $\eta$ Taylor scheme \ref{SWGcom}, respectively.
Assume that $\E\big[|X_0|^q\big]<\infty$ for $q\ge 1$ and
$b^{k},\sigma^{k,j},c^k\in C_P^{\eta+1,2(\eta+1),2(\eta+1)}\big([0,T]\times\R^d\times\R^d;\R\big)$
are Lipschitz continuous for $1\le k\le d$ and $1\le j\le m$.
Let $\tilde{f}_{\alpha}$ with $f(t,x',x)=x$
satisfy
\[
\big|\tilde{f}_{\alpha}^{\beta}(t,x)\big|\le K\big(1+\E[|\beta_t|]+|x|\big),\qquad (t,x)\in[0,T]\times\R^d
\]
for all  $\alpha\in\Gamma_{\eta}\bigcup\mathcal{B}(\Gamma_{\eta})$
with $K>0$ being a constant and $\beta_t$ defined by \eqref{itpppj}.
Then for any function $g\in C_P^{2(\eta+1)}(\R^d;\R)$,
it holds that
\[
\big|\E\big[g(X_T)-g(X_{N})\big]\big|\le C(\Delta t)^\eta,
\]
where $C$ is a positive constant independent of $\Delta t$.
\end{coro}

\section{Numerical examples}\label{sec:numerical}

In this section, we carry out some numerical tests
to verify our theoretical conclusions and
to show the efficiency and the accuracy of the proposed schemes for solving MSDEJs.
For each example, we shall test the Euler scheme \eqref{eul},
the strong order 1.0 Taylor scheme \eqref{mil}, and the weak order 2.0 Taylor scheme \eqref{weak2.0},
respectively.

For simplicity, we adopt the uniform time partition, and the time partition number
$N$ is given by $N=\frac{T}{\Delta t}$.
We denote by
$\mathbb E\big[|X_T-X_{N}|\big]$ and $\big|\mathbb E[X_T-X_{N}]\big|$
the strong errors  and the weak errors  between the exact solution $X_t$ of the MSDEJ \eqref{MSDEJeq}
at time $t=T$ and the numerical solution $X_n$ of the proposed schemes at $n = N$.
The Monte Carlo method is used
to approximate the expectation
$\mathbb E[\cdot]$ appeared in coefficients and errors with sample times $M$.
The ``exact" solution of the MSDEJs is identified with the numerical one
using a small step-size $\Delta t_{\rm exact}= 2^{-12}$.
Moreover, we will test the efficiency of our schemes with respect to the level of
the intensity $\dot{\lambda}$ of the Poisson measure $\mu$ by
the magnitudes of the sample times $M$ and the running time (RT) for different values of $\dot\lambda$.
%

In what follows, we denote by Euler, S-1.0 and W-2.0  the Euler scheme,
the strong order 1.0 Taylor scheme,
and the weak order 2.0 Taylor scheme, respectively.
The convergence rate (CR) with respect to  $\Delta t$ is obtained by using linear
least square fitting to the numerical errors.
In all the tests, we set $T=1.0$.
The unit of RT is the second.


\begin{ex}
Consider the following MSDEJ with $X_0=x_0$:
\begin{equation} \label{ex1-eq}
dX_s=a\big(\E[X_s]+X_s\big)ds+bX_sdW_s+\int_{\msE}ce\big(\E[X_s]+X_{s-}\big)\mu(de,ds),
\end{equation}
where $a$, $b$ and $c$ are constants.
\end{ex}

We set $a=1.25$, $b=0.75$, $c=0.25$, and $X_0=0.1$.
Assume that the jump sizes $\{Y_i,i=1,\dots,N_T\}$
satisfy $Y_i\overset{iid}\sim U(-\frac12,\frac12)$,
which is the uniform distribution on $[-\frac12,\frac12]$.
And we us the Euler scheme \eqref{eul},
the strong order 1.0 Taylor scheme \eqref{mil} and the weak order 2.0 Taylor scheme \eqref{weak2.0}
to solve \eqref{ex1-eq}, respectively.
We have listed the errors and convergence rates of the schemes \eqref{eul},
\eqref{mil} and \eqref{weak2.0} for different intensity $\dot\lambda$ in Tables \ref{table:1a} - \ref{table:1c}, respectively.


\begin{table}[htbp!]\renewcommand{\arraystretch}{1.0}
\centering
\caption{Errors and convergence rates of the Euler scheme}\label{table:1a}
\scalebox{0.9}{\begin{tabular}{|c|ccccc|c|c|c|}
\hline
 \multicolumn{9}{|c|}{Euler\qquad\quad} \\  \hline
${N}$ &  $16$      & $32$      & $64$       & $128$     &  $256$     &  CR    & M   & RT\\ \hline
 $\dot\lambda$ &\multicolumn{8}{c|}{$\E\big[|X(T)-X_{N}|\big]\qquad\qquad~~$}\\  \hline
0.1   & 2.739E-01 & 1.405E-01 & 6.630E-02  & 4.517E-02  & 3.377E-02  & 0.768  & 45  & 2.94\\
0.5   & 1.839E-01 & 9.217E-02 & 4.957E-02  & 3.212E-02  & 2.063E-02  & 0.783  & 65  & 4.18\\
1.0   & 2.174E-01 & 1.236E-01 & 7.031E-02  & 4.236E-02  & 2.479E-02  & 0.781  & 75  & 4.87\\
2.0   & 2.106E-01 & 1.173E-01 & 5.759E-02  & 3.800E-02  & 2.610E-02  & 0.765  & 85  & 5.46\\
3.0   & 1.479E-01 & 8.452E-02 & 5.922E-02  & 3.605E-02  & 2.614E-02  & 0.623  & 100 & 6.14\\ \hline
 \multicolumn{9}{|c|}{$\big|\E[X(T)-X_{N}]\big|\qquad$} \\  \hline
0.1   & 2.192E-01 & 1.222E-01 & 6.122E-02  & 2.696E-02  & 1.034E-02  & 1.099  & 100 & 6.07\\
0.5   & 2.987E-01 & 1.510E-01 & 7.575E-02  & 4.380E-02  & 1.907E-02  & 0.972  & 200 & 12.94\\
1.0   & 3.229E-01 & 1.699E-01 & 8.444E-02  & 4.397E-02  & 1.934E-02  & 1.007  & 300 & 24.51\\
2.0   & 3.401E-01 & 1.880E-01 & 9.122E-02  & 4.380E-02  & 2.007E-02  & 1.027  & 500 & 59.56\\
3.0   & 2.996E-01 & 1.640E-01 & 8.564E-02  & 4.028E-02  & 1.946E-02  & 0.992  & 900 & 107.35\\ \hline
\end{tabular}}
\end{table}

\begin{table}[htbp!]\renewcommand{\arraystretch}{1.0}
\centering
\caption{Errors and convergence rates of the strong order 1.0 Taylor scheme}\label{table:1b}
\scalebox{0.89}{\begin{tabular}{|c|ccccc|c|c|c|}
\hline
 \multicolumn{9}{|c|}{S-1.0\qquad\quad} \\  \hline
${N}$ & $16$      & $32$       & $64$       & $128$      &  $256$     &  CR    & M   & RT\\ \hline
 $\dot\lambda$ &\multicolumn{8}{c|}{$\E\big[|X(T)-X_{N}|\big]\qquad\qquad~~$}\\  \hline
0.1   & 1.679E-01 & 8.939E-02  & 4.560E-02  & 2.292E-02  & 1.126E-02  & 0.976  & 100 & 6.76\\
0.5   & 2.027E-01 & 1.107E-01  & 5.738E-02  & 2.859E-02  & 1.408E-02  & 0.965  & 200 & 13.63\\
1.0   & 2.084E-01 & 1.138E-02  & 5.930E-02  & 2.962E-02  & 1.449E-02  & 0.963  & 400 & 41.21\\
2.0   & 2.138E-01 & 1.164E-01  & 6.046E-02  & 3.053E-02  & 1.477E-02  & 0.964  & 800 & 99.20\\
3.0   & 2.006E-01 & 1.092E-01  & 5.665E-02  & 2.877E-02  & 1.404E-02  & 0.960  & 1000 & 130.32\\\hline
 \multicolumn{9}{|c|}{$\big|\E[X(T)-X_{N}]\big|\qquad$} \\  \hline
0.1   & 1.679E-01 & 8.939E-02  & 4.560E-02  & 2.292E-02  & 1.126E-02  & 0.976  & 100 & 6.76\\
0.5   & 2.027E-01 & 1.107E-01  & 5.738E-02  & 2.859E-02  & 1.408E-02  & 0.965  & 200 & 13.63\\
1.0   & 2.084E-01 & 1.138E-01  & 5.927E-02  & 2.959E-02  & 1.447E-02  & 0.964  & 400 & 41.21\\
2.0   & 2.138E-01 & 1.164E-01  & 6.045E-02  & 3.049E-02  & 1.473E-02  & 0.965  & 800 & 99.20\\
3.0   & 2.006E-01 & 1.092E-01  & 5.657E-02  & 2.867E-02  & 1.394E-02  & 0.962  & 1000 & 130.32\\ \hline
\end{tabular}}
\end{table}

\begin{table}[htbp!]\renewcommand{\arraystretch}{1.0}
\centering
\caption{Errors and convergence rates of the weak order 2.0 Taylor scheme}\label{table:1c}
\scalebox{0.88}{\begin{tabular}{|c|ccccc|c|c|c|}
\hline
 \multicolumn{9}{|c|}{W-2.0\qquad\quad} \\  \hline
${N}$ &  $8$      & $16$      & $32$       & $64$       & $128$      &  CR    & M    & RT\\ \hline
 $\dot\lambda$ &\multicolumn{8}{c|}{$\E\big[|X(T)-X_{N}|\big]\qquad\qquad~~$}\\  \hline
0.1   & 3.686E-02 & 1.329E-02 & 6.709E-03  & 2.672E-03  & 1.271E-03  & 1.203  & 100  & 7.28\\
0.5   & 5.150E-02 & 1.722E-02 & 6.706E-03  & 3.670E-03  & 1.849E-03  & 1.183  & 200  & 16.15\\
1.0   & 4.686E-02 & 1.681E-02 & 6.918E-03  & 3.749E-03  & 2.165E-03  & 1.104  & 500  & 78.36\\
2.0   & 5.022E-02 & 1.974E-02 & 8.703E-03  & 4.663E-03  & 2.825E-03  & 1.039  & 800  & 114.51\\
3.0   & 4.808E-02 & 1.841E-02 & 9.511E-03  & 5.843E-03  & 3.511E-03  & 0.921  & 1500 & 228.86\\ \hline
 \multicolumn{9}{|c|}{$\big|\E[X(T)-X_{N}]\big|\qquad$} \\  \hline
0.1   & 3.438E-02 & 1.013E-02 & 2.618E-03  & 5.500E-04  & 8.881E-05  & 2.140  & 2800 & 450.22\\
0.5   & 4.033E-02 & 1.157E-02 & 2.537E-03  & 7.168E-04  & 1.165E-04  & 2.089  & 3000 & 481.91\\
1.0   & 3.529E-02 & 9.424E-03 & 2.526E-03  & 5.719E-04  & 1.246E-04  & 2.033  & 3500 & 567.14\\
2.0   & 3.945E-02 & 1.101E-02 & 2.587E-03  & 3.550E-04  & 1.099E-04  & 2.193  & 4500 & 826.78\\
3.0   & 3.836E-02 & 1.057E-02 & 2.542E-03  & 6.279E-04  & 4.323E-05  & 2.366  & 5500 & 1201.71\\ \hline
\end{tabular}}
\end{table}
~\\~\\~\\~\\

The numerical results listed in Tables \ref{table:1a} - \ref{table:1c}
show that the Euler scheme \eqref{eul},
the strong order 1.0 Taylor scheme \eqref{mil} and the weak order 2.0 Taylor scheme \eqref{weak2.0}
are stable and accurate for solving the linear  MSDEJ \eqref{ex1-eq}.
Moreover, we can drawn the following conclusions.
\begin{enumerate}
\item The orders of strong convergence of the Euler scheme, the strong order 1.0 Taylor scheme
and the weak order 2.0 Taylor scheme
are 0.5, 1.0 and 1.0, respectively;
\item The orders of weak convergence of the Euler scheme,
the strong order 1.0 Taylor scheme and
the weak order 2.0 Taylor scheme are 1.0, 1.0 and 2.0, respectively;
\item The efficiency of the schemes depends on the level of the intensity $\dot\lambda$
of the Poisson measure $\mu$.
As the intensity $\dot\lambda$ increases, the sample times $M$ and the running time RT increase.
\end{enumerate}
All of the conclusions above are consistent with our theoretical reuslts.

 \begin{ex}
Consider the following general nonlinear MSDEJ
\begin{equation} \label{ex2-eq}\begin{aligned}
dX_t^{0,X_0}=&\; \Big(\big(X_t^{0,X_0}\big)^{5/3}+2\dot\lambda^2\E\big[X_t^{0,x_0}\big]\Big)dt+\frac12\E\big[X_t^{0,x_0}\big]dW_t\\
&\;+\int_{\msE}\frac{e}{2(1+\dot\lambda^2)}\Big(X_t^{0,X_0}+\E\big[\big(X_t^{0,x_0}\big)^2\big]\Big)\mu(de,dt),
\end{aligned}\end{equation}
where $x_0$ and $X_0$ are initial values and $\dot\lambda$ is the intensity of the Poisson measure $\mu$.
\end{ex}

Let the jump sizes satisfy
$Y_i\overset{iid}\sim U(-\frac12,\frac12)$, $i=1,\dots,N_T$.
We use the Euler scheme \eqref{eul},
the strong order 1.0 Taylor scheme \eqref{mil} and the weak order 2.0 Taylor scheme \eqref{weak2.0}
to solve \eqref{ex2-eq}, repsectively.
For simplicity, we set $\dot\lambda=1.0$.
In Tables \ref{table:2a} and \ref{table:2b}, we have listed
the errors and convergence rates of the schemes for different initial values
of $x_0$ and $X_0$.

\begin{table}[htbp!]\renewcommand{\arraystretch}{1.0}
\centering
\caption{Errors and convergence rates of the schemes with $x_0=X_0=0.1$}\label{table:2a}
\scalebox{1.041}{\begin{tabular}{|c|ccccc|c|}
\hline
${N}$     &  $16$     & $32$      & $64$       & $128$     & $256$     & CR        \\ \hline
 \multicolumn{7}{|c|}{$\E\big[|X(T)-X_{N}|\big]$}  \\  \hline
{Euler}   & 3.404E-01 & 1.956E-01 & 1.050E-01  & 5.330E-02 & 3.003E-02 & 0.888  \\
{S-1.0}   & 2.968E-01 & 1.697E-01 & 9.092E-02  & 4.684E-02 & 2.324E-02 & 0.921    \\
{W-2.0}   & 3.207E-02 & 9.929E-03 & 3.726E-03  & 1.848E-03 & 1.110E-03 & 1.213   \\
\hline
 \multicolumn{7}{|c|}{$\big|\E[X(T)-X_{N}]\big|$}  \\  \hline
{Euler}   & 2.802E-01 & 1.605E-01 & 8.578E-02  & 4.349E-02 & 2.110E-02 & 0.935   \\
{S-1.0}   & 2.965E-01 & 1.696E-01 & 9.086E-02  & 4.682E-02 & 2.324E-02 & 0.920    \\
{W-2.0}   & 2.417E-02 & 6.095E-03 & 1.432E-03  & 2.852E-04 & 7.891E-05 & 2.094    \\
\hline
\end{tabular}}
\end{table}

%

\begin{table}[htbp!]\renewcommand{\arraystretch}{1.0}
\centering
\caption{Errors and convergence rates of the schemes with $x_0=0.15$ and $X_0=0.05$}\label{table:2b}
\scalebox{1.041}{\begin{tabular}{|c|ccccc|c|}
\hline
${N}$     &  $16$     & $32$      & $64$       & $128$     & $256$     & CR        \\ \hline
 \multicolumn{7}{|c|}{$\E\big[|X(T)-X_{N}|\big]$}  \\  \hline
{Euler}   & 7.643E-01 & 4.595E-01 & 2.545E-01  & 1.315E-01 & 7.290E-02 & 0.859   \\
{S-1.0}   & 5.989E-01 & 3.559E-01 & 1.945E-01  & 9.997E-02 & 4.898E-02 & 0.906    \\
{W-2.0}   & 8.133E-02 & 3.081E-02 & 1.054E-02  & 4.521E-03 & 2.761E-03 & 1.253    \\
\hline
 \multicolumn{7}{|c|}{$\big|\E[X(T)-X_{N}]\big|$}  \\  \hline
{Euler}   & 5.947E-01 & 3.527E-01 & 1.923E-01  & 9.924E-02 & 4.858E-02 & 0.906   \\
{S-1.0}   & 5.946E-01 & 3.527E-01 & 1.923E-01  & 9.926E-02 & 4.861E-02 & 0.905  \\
{W-2.0}   & 7.884E-02 & 2.870E-02 & 7.801E-03  & 1.535E-03 & 2.047E-04 & 2.140   \\
\hline
\end{tabular}}
\end{table}

From the numerical results in Tables \ref{table:2a} and \ref{table:2b},
we come to the conclusion that the Euler scheme \eqref{eul},
the strong order 1.0 Taylor scheme \eqref{mil} and the weak order 2.0 Taylor scheme
\eqref{weak2.0} are stable and accurate for solving the nonlinear MSDEJ \eqref{ex2-eq}
with different initial values of $x_0$ and $X_0$.
Tables \ref{table:2a} and \ref{table:2b} also show that the orders of strong convergence
of the schemes \eqref{eul}, \eqref{mil} and \eqref{weak2.0}
are 0.5, 1.0 and 1.0, respectively, and the orders of weak convergence
are 1.0, 1.0 and 2.0, respectively,
which verify again our theoretical conclusions.
~\\

Note that the efficiency of the proposed schemes, especially the high order ones depends
on the level of the intensity of the Poisson measure.
This is mainly due to the existence of the double integrals
involving the Poisson measure
$I_{\{-1,-1\}}$, $I_{\{1,-1\}}$, $I_{\{-1,1\}}$, $I_{\{0,-1\}}$ and $I_{\{-1,0\}}$,
the computation complexity of which is dependent on the number of jumps.
Hence, to construct more efficient high order schemes for MSDEJs, we will focus on
the jump-adapted methods in our future work,
which avoid the integrals involving the Poisson measure.
\section{Conclusions}\label{sec:con}

In this paper, we developed the It\^o formula and
the It\^o-Taylor expansion for MSDEJs,
then based on which we proposed the strong order $\gamma$
and the weak order $\eta$ It\^o-Taylor schemes for solving MSDEJs.
We rigorously proved
the error estimates of the proposed schemes, which show that
the order of strong convergence of
the strong order $\gamma$ Taylor scheme and
the order of weak convergence of the weak order $\eta$ Taylor scheme
are $\gamma$ and $\eta$, respectively.
Numerical experiments verify our
theoretical conclusions and indicate that the efficiency of
the schemes depends on the level of
the intensity of the Poisson measure.
In the future work, we shall consider the jump-adapted methods for MSDEJs.

\end{document}